\documentclass[10pt,amssymb]{article} 
\usepackage{amsfonts,latexsym}

\newcommand{\R}{\mathbb R}  
\newcommand{\Z}{\mathbb Z}
\newcommand{\N}{\mathbb N}  

\newcommand{\LL}{\mathbb L}

\newcommand{\ms}{\medskip}

\newcommand{\del}{\partial}  
\newcommand{\CI}{{\mathcal C}^\infty}
\newcommand{\e}{\varepsilon}  

\newcommand{\bx}{{\bf x}} 

\newcommand{\Si}{\Sigma}

\def\ds{\displaystyle}

\newtheorem{theorem}{Theorem}  
\newtheorem{lemma}{Lemma}
\newtheorem{proposition}{Proposition}  
\newtheorem{corollary}{Corollary}
\newtheorem{definition}{Definition}  
\newtheorem{remark}{Remark}
 
\begin{document} 
 
\title{Connected sums of constant mean curvature surfaces in  Euclidean 
3 space.}

\author{Rafe Mazzeo\thanks{Supported by 
the NSF under Grant DMS-9626382} \\ Stanford University \and
 Frank Pacard\thanks{Supported by AIM (American Institute of  
Mathematics)} \\ Universit\'e Paris XII  \and 
Daniel Pollack\thanks{Supported by the NSF under Grant DMS-9704515}  
\\ University of  Washington}

\date{May 11, 1999}

\maketitle

\section{Introduction and statement of the results}
\label{se:1}

Amongst the recent developments in the study of embedded complete
minimal and constant mean curvature surfaces in $\R^3$ is the
realization that these objects are far more robust and flexible than is
apparent from  their Weierstrass representations. Our aim in this paper
is to prove a `gluing theorem', which states roughly that if two
(appropriate) constant mean curvature surfaces are juxtaposed, so that
their tangent planes are parallel and very close to one another, but
oppositely oriented, then there is a new constant mean curvature surface 
quite near to this configuration (in the Hausdorff topology), 
but which is a topological connected 
sum of the two surfaces. We shall explain what we  mean by appropriate, 
or at least give our preliminary interpretation of it, in the next
paragraph. Throughout this paper, the acronym CMC shall mean a surface
with constant mean curvature equal to one (or minus one depending on the
orientation).

The simplest context for our result is when we are given two orientable,
immersed, compact CMC surfaces, $\Sigma_1$ and $\Sigma_2$, with nonempty
boundary. Suppose that we have applied a rigid motion to each of these
surfaces so that $0 \in \Sigma_1 \cap \Sigma_2$ and $T_0 \Sigma_1 = T_0
\Sigma_2$ is the $x\, y$-plane. (These surfaces may intersect elsewhere, but 
that is irrelevant for our considerations.)  We now define the 
orientation on these surfaces so that at $0$ the oriented 
unit normal $\nu_1$ of $\Sigma_1$ equals $(0,0,1)$, while the oriented 
unit normal $\nu_2$ of $\Sigma_2$ equals $(0,0,-1)$. Let us assume that
with this orientation the two surfaces have the same mean curvature $H_0$
(so either $H_0=1$ or $H_0=-1$ for both of the surfaces). 
We shall prove that there is a `geometric connected sum' of
these two surfaces, which may be thought of as a desingularization of
this configuration. Moreover, the boundary of this desingularization
will be the union of the boundaries of the $\Sigma_i$, each possibly
transformed by a small rigid motion.

In order to state this first result rigorously, we make the following
definition:
\begin{definition}
A compact CMC surface $\Sigma$ with boundary is said to  be {\bf
nondegenerate} if there are no Jacobi fields on $\Sigma$ which vanish on
$\del \Sigma$.  Namely, if $w : \Sigma \longrightarrow \R$ is a ${\cal
C}^{2, \alpha}$  solution of
\[
\Delta_\Sigma w + |{\bf A}_\Sigma|^2 w =0, \qquad \left. w \right|_{\del
\Sigma} = 0,
\]
then $w=0$. Here ${\bf A}_\Sigma$ is the second fundamental form of
$\Sigma$.
\label{de:1}
\end{definition}

\begin{theorem} ({\bf Connected sum theorem})
Let $\Sigma_1$ and $\Sigma_2$ be two compact, smooth, immersed,
orientable, {\bf nondegenerate} CMC surfaces with boundary. Assume that
these surfaces are positioned and oriented as above and have the same mean 
curvature $H_0$.  Then there exist
an $\e_0>0$ and a one-parameter family  of surfaces $S_\e$, for $\e \in
(0, \e_0]$, satisfying the following  properties:
\begin{enumerate}
\item $S_\e$ is a smooth, immersed CMC surface with boundary.

\item There are rigid motions $\tau_1$ and $\tau_2$ of $\R^3$, depending
on $\e$, such that $\partial S_\e = \tau_1(\partial\Sigma_1) \cup
\tau_2(\partial\Sigma_2)$.

\item For any fixed $R>0$, the surface $S_\e \cap \, [{\R}^3 \setminus
B_R]$  converges in the $\CI$ topology to $[\Sigma_1 \cup \Sigma_2 ] \,
\cap \, [{\R}^3 \setminus B_R]$ and  $\del S_\e$ converges in the  $\CI$
topology to $\del\Sigma_1\cup\del\Sigma_2$.

\item The dilated surface $\e^{-1}S_\e$ converges in the $\CI$
topology on any compact set to a catenoid with vertical axis.
\end{enumerate}
\label{th:1}
\end{theorem}

\begin{remark}
There are actually two geometrically distinct families of
surfaces $S_\e$ which are constructed here, corresponding
to the two choices $H_0 = \pm 1$. To better visualize these,
consider two small spherical caps intersecting at the origin
and both tangent to the $xy$ plane.  Assume that these surfaces
are oriented oppositely to one another so that one is below the $xy$ 
plane and the other is above. If their normals are pointing outward 
(so that their mean curvatures are both $-1$), then the new 
surfaces $S_\e$ are embedded and very much resemble a neighbourhood
of the neck region in an embedded Delaunay surface (an unduloid).
Of course, by reversing the orientation of these resulting surfaces we 
obtain surfaces with  mean curvature $=+1$.
On the other hand, if the initial orientations are reversed so that
the mean curvatures are both $+1$, then the resulting $S_\e$ are
only immersed, and resemble the neck regions in the immersed Delaunay
surfaces of nodoid type. 
\end{remark}

We also obtain additional geometric information about the surfaces
$S_\e$, in particular that their geometry is well-controlled as $\e \to
0$.
\begin{proposition} ({\bf Embeddedness})
Under the assumptions of the previous theorem, assume further that
$[\Sigma_1 \cup \Sigma_2 ]\setminus \{0\}$ is embedded. Then for
one of the two choices of $H_0$, and for $\e$ sufficiently small, 
the surface $S_\e$ is embedded.
\label{pr:1}
\end{proposition}
We also obtain estimates on the rate of convergence of $S_\e$ to
$\Sigma_1\cup \Sigma_2$.
\begin{proposition} ({\bf Distance from $\Sigma_1\cup \Sigma_2$ to $S_\e$})
Again under the assumptions of the previous theorem, there exists a
constant $c>0$ such that
\[
\mbox{\em dist} \, (S_\e , \Sigma_1 \cup \Sigma_2) \equiv \max \left(
\sup_{p \in \Sigma_1 \cup \Sigma_2}\mbox{\em dist}\, (p, S_\e),
\sup_{q \in S_\e} \mbox{\em dist}\, (q, \Si_1 \cup \Si_2) \right)
\leq c \, \e \, |\log \e|.
\]
\end{proposition}

There are various ways to generalize these results. First assume  that
$\Sigma_1$ and $\Sigma_2$ are two smooth, oriented CMC surfaces with
boundary, which are nondegenerate. If $p_i \in \Sigma_i$,  $i=1,2$, then
we apply a rigid motion to each surface so  that $p_1=p_2=0$ and
$T_0\Sigma_i$ is the $x \, y$-plane (with  opposite orientations as
above, and same mean curvature for the two surfaces).
 Next rotate the surface $\Sigma_2$  about the $z$-axis by an
angle $\theta \in S^1$. This gives a five parameter family of initial
configurations $\Sigma_1 \sqcup  \Sigma_2 (p_1, p_2, \theta)$. The precise
definition of  $\Sigma_1 \sqcup \Sigma_2 (p_1, p_2, \theta)$ will be given
in \S 16. Applying Theorem~\ref{th:1} to desingularize each of these
configurations adds an additional parameter, and we obtain the six
parameter family $S_\e (p_1, p_2, \theta)$.

It turns out that this family depends smoothly on all six parameters.
We will not prove this explicitly in this paper, in order to  keep the
technicalities to a minimum; however, the proof is not hard to deduce
from our arguments. This dimension count is closely related to the
question of whether the solutions $S_\e(p_1, p_2, \theta)$ are
nondegenerate. For if this is the case for one of these surfaces, then
the implicit function theorem gives a six dimensional smooth family of
CMC surfaces in a neighbourhood of that surface. Unfortunately we can
only prove that these surfaces  are nondegenerate for generic choices of
parameters.
\begin{proposition} ({\bf Generic nondegeneracy property})
There is a (singular) codimension one analytic set ${\cal S}$ in
$\Si_1 \times \Si_2 \times S^1$ such that the surface 
$S_\e(p_1,p_2,\theta)$ is nondegenerate provided $\e$ is
sufficiently small and $(p_1,p_2,\e) \notin {\cal S}$. 
The set ${\cal S}$ is the union of the locus of points
satisfying a quadratic polynomial equation in $\cos \theta$
with coefficients depending on the principal curvatures of
the surfaces at $p_i$, together with a set ${\cal C} \times S^1$,
where ${\cal C}$ is the product of the locus of points on
the two surfaces where $\Si_1$ is umbilic (and hence the principal curvatures
are equal to $(H_0/2,H_0/2)$) and the principal curvatures on $\Si_2$ are equal 
to $(-H_0/2,3 H_0/2)$, or vice versa. ${\cal C} \times S^1$ has dimension
less than four unless $\Si_1$ or $\Si_2$ is a subdomain of the sphere.
\label{pr:2}
\end{proposition}
Notice that we have defined the mean curvature to be the sum of the 
two principal curvatures, not the average.

One important application of this result is that if $S_\e (p_1, p_2,
\theta)$ is nondegenerate, then one can use it as one of the `summands'
in another application of Theorem~\ref{th:1}, and so the connected sum
procedure may be iterated. Thus, for example since certain subdomains
of the sphere or cylinder with nonempty boundary are 
nondegenerate, we may glue together arbitrarily many copies of them. 

Gluing constructions for geometric objects are by now well understood
and even somewhat commonplace, and they have been used to solve a number
of diverse problems. Even in the context of CMC and minimal surfaces,
there are many results. The pioneering work in this area was that
of N. Kapouleas, cf. \cite{K1}, \cite{K2},  \cite{K3} and \cite{K4}.
Recently, S.D. Yang \cite{Y} has proved a connected sum theorem for 
complete minimal surfaces of finite total curvature. The methods
here could equally well be used to prove that result (or indeed, his
methods could be used in the present context), but although 
Yang requires nondegeneracy of his minimal summands, he does not
discuss the question of nondegeneracy of the final surface at 
all, and it is not clear how it could be obtained by that
approach. The issue of nondegeneracy is quite important in the 
moduli space theory, cf. \cite{KMP} and the recent work \cite{GKS}. 

The results in this paper are also close in spirit to the connected sum
theorem in the scalar curvature context in \cite{MPU}, but the methods
there are much simpler.  The method of proof here is inspired by the
recent work of the first and second authors \cite{MPa} on the
construction of CMC surfaces with finitely many Delaunay ends. We now
briefly comment on our construction, pointing out its novel features.

The usual steps in such a construction would be to first build a family
of approximate solutions, depending on a parameter $\e > 0$.
These approximate CMC surfaces would consist of the surfaces $\Si_1$ 
and $\Si_2$ and a catenoidal neck, joined together with cutoff functions,
and would  converge to the singular configuration $\Sigma_1 \cup \Sigma_2$ 
as $\e \to 0$. They would then be perturbed, when $\e$ is
sufficiently small, to obtain the desingularized CMC surface. This step
involves a careful analysis of the Jacobi operator of these approximate
solutions, uniformly as $\e \to 0$.

We proceed somewhat differently here. Our building blocks are
the same, namely the surfaces $\Si_1$ and $\Si_2$ and a small
`neck region' of a catenoid. Roughly speaking, we construct
perturbations of each of these components which are themselves CMC
surfaces with boundary in such a way that the
Cauchy data matches across the boundary. The boundary here
consists of the small curves produced by excising small balls
around the points $p_1$ and $p_2$ as well as the boundaries
of the truncated catenoid. A very important point 
is that we first perturb each of the surfaces $\Sigma_i$ by
adding in the normal direction $\e$ times the Green function 
for the Jacobi operator with pole at $p_i$. This is the precise
point where we use nondegeneracy of the surfaces $\Si_i$,
and has the important effect of making the local geometry of 
$\Sigma_i$ near $p_i$ insignificant. The catenoid (scaled
by $\e$) and these surfaces are then truncated at just the
right scale so that their boundaries fit together as well
as possible. 

In the main step of the construction, we construct the infinite
dimensional families of CMC surfaces which are normal graphs over each
of these component pieces. This is done by a simple contraction
mapping argument. As already intimated, we
analyze the Cauchy data of the surfaces in these infinite dimensional 
families at the boundary curves arising from the truncations.
We show using degree theory that this Cauchy data may be
matched, and hence that the desired CMC surface may be
constructed. A substantial advantage of this method is that no 
extraneous cutoff functions are introduced. Because of the high degree of 
nonlinearity of the problem, these are typically the cause of many 
technical complications.

A more detailed guide to the contents is as follows. In \S\S 3, 4 and 5 
we define the truncations of the rescaled catenoid, study the Jacobi 
operator around these surfaces and construct the family of nearby CMC 
surfaces, respectively. \S 6 collects some facts about the
mean curvature operator for graphs and in \S 7 we discuss the
perturbation of the surfaces $\Si_i$ by their Green functions.
\S 8 contains some technical facts about some geometric
modifications of these surfaces arising (mostly) from rigid motions.
Then in \S \S 9 and 10 we study the Jacobi operator on these
modified surfaces and then construct the family of nearby
CMC surfaces. The Cauchy data maps for each of these
components are discussed at the end of \S \S 5 and 10. In the brief
\S 11 we adapt the previous results to our specific needs.
Finally, the degree theory argument for matching the Cauchy data is given
in \S 12.  The remaining sections, \S \S 13, 14 and 15, are devoted
to the analysis of nondegeneracy. \S 13 contains some technical
facts which are needed later, certain Jacobi fields on $S_\e$ 
are discussed in \S 14, and using this nondegeneracy is proved 
in \S 15. There are also three brief appendices containing
various analytic facts which are required in various places
throughout the paper.

The techniques developed here apply immediately to establish a general 
connected sum theorem for complete, noncompact, embedded CMC
surfaces.  For such surfaces there is a natural notion of nondegeneracy 
which in particular 
follows from the nonexistence of square integrable Jacobi fields 
(see \cite{KMP}, \cite{MPo} and \cite{MPa}).  Examples of such surfaces 
are given by the classical Delaunay surfaces \cite{D} which are CMC 
surfaces of revolution, and also the surfaces constructed more recently 
in \cite{MPa}. 
   
In particular, this theorem allows us  to glue together any two embedded 
Delaunay surfaces to produce new embedded four-ended CMC surfaces.
As in \cite{MPU} the resulting surfaces will be asymptotic to the
original Delaunay surface on one end of each pair, the other 
being asymptotic to a small perturbation of the corresponding  
end (here ``small'' is understood within the 6 dimensional family
of Delaunay surfaces which includes those generated by rigid motions).  
Since, in this context, nondegeneracy of the resulting surfaces also
holds generically, the process may be iterated. This produces 
families of complete CMC surfaces which are quite different from 
the previously known examples.  
Moreover, this nondegeneracy together with the control
on the free parameters in our construction allows us to produce an open
subset in the moduli space of complete embedded surfaces with $2k$-ends.  
This open set is actually a collar neighbourhood of  certain boundary 
components in the moduli space.  
Precise statements of
these results along with applications to the study of the moduli space 
itself are given in \cite{MPaP}.

\section{Notation}
\label{se:2}

In this brief section we record some notation that will be used
frequently, throughout the rest of the paper, and without comment.
First, $\lambda : {\R} \longrightarrow [0,1]$ will denote a smooth
cutoff function satisfying
\begin{equation}
\lambda \equiv 1 \qquad \mbox{if} \qquad t  > 1 \qquad \mbox{and}\qquad
\lambda \equiv 0 \qquad \mbox{if} \qquad t  < 0.
\label{eq:0.1}
\end{equation}

\noindent
Next, if $\ds \phi =\sum_{n \in {\Z} } a_n \, e^{in \theta} \in H^{1}
(S^1)$, then we define
\begin{equation}
|D_\theta| \phi \equiv  \sum_{n \in {\Z} } |n| \, a_n \, e^{in \theta}
\in L^2 (S^1).
\label{eq:0.2}
\end{equation}
This is, of course, just $\sqrt{-\Delta} \,\ds \phi$.

\noindent
Finally we define orthogonal projections $\pi'$ and $\pi''$ on
$L^2(S^1)$ as follows: for \[ \phi(\theta) =\sum_{n \in {\Z}} a_n e^{in
\theta} \in L^{2} (S^1),
\]
we set
\begin{equation}
\pi'(\phi) = \sum_{|n| \geq 1} a_n e^{in \theta}, \qquad \mbox{\rm and}  \qquad   \pi''(\phi) =  \sum_{|n| \geq 2} a_n e^{in \theta} \in L^2 (S^1).
\label{eq:0.3}
\end{equation}

\section{The rescaled catenoids}
\label{se:5}

The standard catenoid $\Sigma^c$ has the following standard
parametrization
\begin{equation}
{\bx}^c(s, \theta) = (\cosh s \cos \theta, \cosh s \sin \theta , s),
\qquad (s,\theta) \in {\R}\times S^1.
\label{eq:3.1}
\end{equation}
$\Sigma^c$ may be divided into two pieces, denoted $\Sigma^c_{\pm}$,
which are defined to be the image by ${\bx}^c$ of $({\R}^{\pm} \times
S^1)$, respectively. We may also parameterize the lower half
$\Sigma^c_-$ by
\begin{equation}
{\R}^2 \setminus B_1 \ni (x,y) \longrightarrow (x,y, -\log r - \log 2  +
{\cal O}(r^{-2})) \quad \mbox{\rm as}\ \ r \to \infty.
\label{eq:3.2}
\end{equation}
Here, as usual, $r = (x^2+y^2)^{1/2}$. For any $\e>0$, we define the
rescaled catenoid $\Sigma_\e^c$ by scaling $\Sigma^c$ by the factor
$\e$ and translating by $-\e \, \log \e + \e \,\log 2$ along the  
$z$-axis. $\Sigma_\e^c$ is parameterized by
\begin{equation}
{\bx}_\e^c (s, \theta) = (\e\,\cosh s\,\cos\theta,\e\,\cosh s\,\sin
\theta ,\e\, s - \e \, \log \e + \e \, \log 2), (s,\theta) \in
{\R}\times S^1.
\label{eq:3.3}
\end{equation}
Again $\Sigma_\e^c$ decomposes into two pieces, $\Sigma^c_{\e,\pm}$.  By
(\ref{eq:3.2}) we may parametrize $\Sigma^c_{\e,-}$ either as
\[
{\R}^2 \setminus B_1 \ni (x,y)  \longrightarrow (\e\, x,\e\, y,  - \e\,
\log r - \e \, \log \e + {\cal O}(\e \,r^{-2})),
\]
or equivalently (replacing $(\e x, \e y)$ by $(x,y)$),
\begin{equation}
{\R}^2 \setminus B_\e \ni (x,y) \longrightarrow (x,y, - \e \log r  +
{\cal O}( \e^{3} \, r^{-2} )).
\label{eq:3.4}
\end{equation}
The simplicity of this final parametrization is why we  introduced the
translation along the $z$-axis in the first place.

Finally, consider all surfaces near to the rescaled catenoid
$\Sigma_\e^c$ which may be written as normal graphs off of it.  Since
the outer unit normal of $\Sigma_\e^c$ is given by
\begin{equation}
n (s, \theta) = \frac{1}{\cosh s} (\cos \theta, \sin \theta, - \sinh s ),
\label{eq:5.2}
\end{equation}
each of these surfaces may parameterized as
\begin{equation}
(s, \theta) \longrightarrow {\bx}_\e^c (s, \theta) + w(s,\theta)\,  n(s,
\theta),
\label{eq:5.1}
\end{equation}
for some function $w \in {\mathcal C}^2(\R^- \times S^1)$, which is
suitably small. We prove in Appendix II that the linearized mean
curvature operator about $\Sigma_\e^c$, i.e. at $w=0$, is given by
$-(\e \cosh s)^{-2} {\cal L}$, where
\begin{equation}
{\cal L}  w \equiv   \del^2_{ss} \, w + \del^2_{\theta\theta}\, w   +
\frac{2}{\cosh^2 s}\, w.
\label{eq:5.3}
\end{equation}
As usual, we call this the Jacobi operator. We also prove in this
appendix that the mean curvature of the surface  parameterized by
(\ref{eq:5.1}) is given by an expression of the form
\begin{equation}
\begin{array}{rlll}
H_w = - \ds \frac{1}{\e^2 \cosh^2 s}\, {\cal L} w & + & \ds \frac{1} {\e
\cosh^2 s}  Q_\e'\left(\frac{w}{\e \cosh s},\frac{\nabla w}{\e \cosh
s},\frac{\nabla^2  w}{\e \cosh s}\right) \\[3mm] & + & \ds \frac{1}{\e
\cosh s} Q_\e''\left(\frac{w}{\e \cosh s}, \frac{\nabla w}{\e \cosh
s},\frac{\nabla^2 w}{\e \cosh s}\right),
\end{array}
\label{eq:5.44}
\end{equation}
where $Q'_\e$ and $Q''_\e$ are functions which are bounded in
${\mathcal C}^k ([-s_\e, s_\e] \times S^1)$ for all $k$,  uniformly in
$\e$. These functions also satisfy
\begin{equation}
Q'_\e (0,0,0)= Q''_\e (0,0,0) = 0 \qquad \mbox{and}\qquad  \nabla Q'_\e
(0,0,0)= \nabla Q''_\e (0,0,0) =0 ,
\label{eq:5.55}
\end{equation}
and in addition
\begin{equation}
\nabla^2 Q''_\e (0,0,0)=0.
\label{eq:5.66}
\end{equation}

\section{The mean curvature operator linearized about 
the truncated catenoid}
\label{se:7}

In the next section we shall study the space of CMC surfaces in a
neighbourhood of the truncated catenoid. This analysis depends on a good
understanding of the linearization of the mean curvature operator, or
equivalently, of the operator ${\cal L}$ of (\ref{eq:5.3}), on
arbitrarily large truncations of the catenoid. We consider this now.

The mapping properties of ${\cal L}$ are best stated in terms of the
following weighted spaces:

\begin{definition} 
For each $k \in {\N}$ and  $\alpha \in (0,1)$, let $|w|_{k,\alpha, [s, s
+1]}$, denote the usual ${\cal C}^{k,\alpha}$ H\"older norm on the  set
$[s,s+1] \times S^1$. Then for any $\delta\in {\R}$,
\[
{\cal C}^{k,\alpha}_\delta({\R} \times S^1)  =  \left\{ w \in  {\cal
C}^{k,\alpha}_{\mathrm loc}(\R \times S^1):
||w||_{k,\alpha,\delta}\equiv \sup_{s \in {\R}}\,\left[ (\cosh
s)^{-\delta}  |w|_{k,\alpha,[s,s+1]} \right] < \infty \right\}.
\] 
\label{de:2}
\end{definition} 
For any closed interval $I\subset\R$, we denote the restriction  of
${\cal C}^{k,\alpha}_\delta({\R} \times S^1) $ to $I\times S^1$ by
${\cal C}^{k,\alpha}_\delta(I \times S^1)$, endowed with the induced
norm.

\begin{proposition}
Fix $\delta \in (1,2)$. Then for any $s_0 \in {\R}^+$ there  exists an
operator
\[
{\cal G}_{s_0} : {\cal C}^{0,\alpha}_{\delta}([-s_0, s_0]\times S^1)
\longrightarrow  {\mathcal C}^{2,\alpha}_{\delta}([-s_0, s_0] \times S^1)
\]
such that for any $f \in {\cal C}^{0,\alpha}_{\delta}  ( [-s_0, s_0]
\times S^1)$, the function $w={\cal G}_{s_0}(f)$ solves
\begin{equation} 
\left\{
\begin{array}{rllll} 
{\cal L} w & = & \ds  f \qquad & \mbox{in} \quad (- s_0,s_0) \times
S^1\\[3mm]  \pi'' w & =  &  0 \qquad & \mbox{on}\quad \{\pm s_0\} \times
S^1 .
\end{array} 
\right.
\label{eq:5.4}
\end{equation} 
Moreover, $||{\cal G}_{s_0}(f)||_{2, \alpha, \delta} \leq c \, ||f||_{0,
\alpha, \delta },$ for some constant $c>0$ independent of $s_0$.
\label{pr:6}
\end{proposition}

\begin{remark}
The right inverse ${\cal G}_{s_0}$ with these properties is not uniquely
defined. We shall always use the one constructed in the proof below.
\end{remark}

\smallskip

\noindent
{\bf Proof:} Assume that $|f(s,\theta)| \leq (\cosh s)^{\delta}$.  Now
decompose both $w$ and $f$ into Fourier series
\[
w= \sum_{n \in {\Z}} w_n (s) e^{in \theta} \qquad \mbox{and} \qquad  f=
\sum_{n \in {\Z}} f_n (s)  e^{in \theta}.
\]
For $|n|\geq 2$, $w_n$ must solve
\[
\ddot{w}_n -n^2 w_n + \frac{2}{\cosh^2 s} w_n = f_n \qquad \mbox{in}
\qquad  |s| < s_0, \qquad w_n (\pm s_0)=0.
\]
The dots represent differentiation with respect to $s$.

Since $|n|\geq 2$,
\[
L_n = \frac{d^2\,}{ds^2} -n^2  + \frac{2}{\cosh^2 s}
\]
satisfies the maximum principle, so that if $w$ is defined on some
interval  $[s_1,s_2]\subset {\R}$ and if $w(s_1) \geq 0$, $w(s_2) \geq
0$ and  $L_n w \leq 0$ on  $(s_1,s_2)$, then $w \geq 0$ in $[s_1,s_2]$.
We obtain the solution of $L_n w_n = f_n$ by the method of sub- and
supersolutions once we have constructed an appropriate barrier
function. But
\[
L_n (\cosh s)^{\delta} = \left(  (\delta^2 - n^2) \cosh^2s + 2+\delta -
\delta^2 \right) (\cosh s)^{\delta-2},
\]
and then, since $\delta \in (1,2)$,
\[ 
(\delta^2 - n^2)\cosh^2s + 2+\delta -\delta^2 \leq  - (n^2 -2 - \delta)
\cosh^2s.
\]
Therefore, since $|f_n(s)| \leq (\cosh s)^{\delta}$, we have that
\[
L_n(w_n - (n^2 -2- \delta)^{-1}\cosh^\delta s) \geq 0
\]
\[
L_n(w_n + (n^2 -2- \delta)^{-1}\cosh^\delta s) \leq 0.
\]
We conclude that the solution $w_n$ exists and satisfies
\begin{equation}
|w_n(s)| \leq \frac{1}{n^2-2-\delta} (\cosh s)^{\delta}.
\label{eq:5.5}
\end{equation}

Next we obtain the solution and estimates when $n=0, \pm 1$.  This is
straightforward since we know homogeneous solutions of $L_n$ explicitly
for these values of $n$. In fact, $L_0 \tanh s =0$ and $L_{\pm 1}
(\cosh s )^{-1} = 0$. Therefore, by `variation of constants', we obtain
the solutions
\begin{equation}
w_{0}(s) = \tanh s \int_0^s \tanh^{-2} t \int_0^t \tanh u \, f_{0}(u)\,
du\, dt,
\label{eq:5.6}
\end{equation}
and
\begin{equation}
w_{\pm 1} (s) = \cosh^{-1} s \int_0^s \cosh^2 t \int_0^t \cosh^{-1} u \,
f_{\pm 1}(u)\, du\, dt.
\label{eq:5.7}
\end{equation}
Straightforward estimates using these formul\ae\ and the fact that $|f_n
(s)|\leq (\cosh s)^{\delta}$, $n = 0, \pm 1$, gives
\begin{equation}
|w_0(s)| + |w_{\pm 1}(s)| \leq c \, (\cosh s)^{\delta} ,
\label{eq:5.8}
\end{equation}
for some constant $c>0$ independent of $s_0$.

To finish the proof we must amalgamate these estimates. But the
coefficient on the right in (\ref{eq:5.5}) is summable in $n$, and so we
easily see that  $|w(s,\theta)|\leq c\, (\cosh s)^{\delta}$. The
estimates for the  derivatives of $w$ are then obtained by Schauder
theory.  \hfill $\Box$

\medskip

Using a similar technique, we prove the
\begin{proposition}
For each $s_0 >0$ there exists an operator
\[
{\cal P}^0_{s_0}:\left(\pi''\left({\cal C}^{2,\alpha}(S^1)\right)
\right)^2 \longrightarrow {\mathcal C}^{2,\alpha}_{2}([-s_0, s_0]
\times S^1)
\]
such that for all $\phi''_\pm \in \pi'' \left( {\cal C}^{2,\alpha}(S^1)
\right)$, the function $w= {\cal P}^0_{s_0}(\phi''_+, \phi''_-)$ solves
\begin{equation} 
\left\{
\begin{array}{rllll} 
\Delta w & =  & 0 \qquad & \mbox{\rm in}\quad (- s_0,s_0) \times
         S^1\\[3mm]  w & =  &  \phi''_\pm \qquad & \mbox{\rm on} \quad
         \{ \pm s_0 \}\times S^1.
\end{array}  
\right.
\label{eq:5.1111}
\end{equation} 
We also have $||{\cal P}^0_{s_0}(\phi_+,\phi_-)||_{2,\alpha,2} \leq c \,
e^{-2 s_0} \, ( ||\phi''_+ ||_{2,\alpha} +||\phi''_- ||_{2,\alpha})$ for
some $c>0$ independent of $s_0$.
\label{pr:88}
\end{proposition}
{\bf Proof:} By linearity, we may assume that $||\phi''_+ ||_{2,\alpha}
+ ||\phi''_- ||_{2,\alpha} \leq 1$. Again, we decompose $w$ into Fourier
series
\[
w= \sum_{|n|\geq 2} w_n (s) \, e^{in \theta} ,
\]
and obtain the solution by the method of sub- and supersolutions once we
have  constructed an appropriate barrier function. But
\[
\Delta \left( (\cosh s)^n e^{in\theta} \right) = - n (n-1) \, 
(\cosh s)^{n-2}\leq 0.
\]
Therefore, $s\rightarrow (\cosh s_0)^{-n} (\cosh s)^n $ can be used as a
barrier function.  We conclude that the solution $w_n$ exists and
satisfies
\begin{equation}
|w_n(s)| \leq (\cosh s_0)^{-n}\, (\cosh s)^n.
\label{eq:5.555}
\end{equation}
From this it is easy to get the estimate
\[
|w(s,\theta)|\leq c \, (\cosh s_0)^{-2} \, (\cosh s)^2
\]
for all $|s| \leq s_0 - 1$.  The rest of the proof is now obvious and
left to the reader. \hfill $\Box$

\medskip

Using the previous results, we also get
\begin{proposition}
Fix $\delta \in (1,2)$. Then for each $s_0 >0$ there exists an operator
\[
{\cal P}_{s_0}:\left(\pi''\left({\cal C}^{2,\alpha}(S^1)\right)
\right)^2 \longrightarrow {\mathcal C}^{2,\alpha}_{\delta}([-s_0, s_0]
\times S^1)
\]
such that for all $\phi''_\pm \in \pi'' \left( {\cal C}^{2,\alpha}(S^1)
\right)$, the function $w= {\cal P}_{s_0}(\phi''_+, \phi''_-)$ solves
\begin{equation} 
\left\{
\begin{array}{rllll} 
{\cal L} w & =  & 0 \qquad & \mbox{\rm in}\quad (- s_0,s_0) \times
         S^1\\[3mm]  w & =  &  \phi''_\pm \qquad & \mbox{\rm on} \quad
         \{ \pm s_0 \}\times S^1.
\end{array}  
\right.
\label{eq:5.11}
\end{equation} 
We also have $||({\cal P}_{s_0} -{\cal
P}^0_{s_0})(\phi_+,\phi_-)||_{2,\alpha,\delta} \leq c \, e^{-2 s_0} \, (
||\phi''_+ ||_{2,\alpha} +||\phi''_- ||_{2,\alpha})$ for some $c>0$
independent of $s_0$.
\label{pr:8}
\end{proposition}

{\bf Proof:} Set $w = w_0 + v$ where the functions $w_0$ and $v$
are given by $w_0 = {\cal P}^0_{s_0}(\phi''_+, 
\phi''_-)$ and $v = -{\cal G}_{s_0}{\cal L} w_0 = - 2{\cal G}_{s_0} 
\left( (\cosh s)^{-2} w_0 \right)$, $v(\pm s_0,\theta) = 0$. 
Then the estimate $||w_0||_{2,\alpha,2} \leq c \, e^{-2 s_0}\,
(||\phi''_+||_{2,\alpha} + ||\phi''_-||_{2,\alpha})$ and an
application of Proposition~\ref{pr:6} give the estimate for $v$ and 
finishes the proof.\hfill $\Box$

\ms

To simply notation we shall henceforth write ${\cal
P}^0_{s_0}(\phi''_{\pm})$, ${\cal P}_{s_0}(\phi''_{\pm})$ and
$||\phi''_{\pm}||_{2, \alpha}$ in place of the longer versions above.

\section{CMC surfaces near the truncated catenoid}
\label{se:8}

Now and hereafter, we set
\begin{equation}
s_\e = -\frac{1}{4}\log \e .
\label{eq:6.1}
\end{equation}
Use the parametrization (\ref{eq:3.3}) for the rescaled catenoid. Its
outer unit normal $n(s,\theta)$ at ${\bf x}_\e^c (s,\theta)$ is then
given  by (\ref{eq:5.2}). Define a smooth, strictly monotone function
$\xi_\e : {\R} \longrightarrow [-1,1]$ by
\begin{equation}
\xi_\e (s)=-\left(1-\lambda (s_\e-1-|s|)\right)\, \frac{s}{|s|} -
\lambda (s_\e-1-|s|)\,\tanh s,
\label{eq:6.2}
\end{equation}
Thus $\xi_\e (s) = -\ds \frac{s}{|s|}$ for $|s| \geq s_\e - 1$ and
$\xi_\e (s) = - \tanh s$ for $ |s| \leq s_\e -2$. Now consider the
vector field
\begin{equation}
\bar{n}_\e (s, \theta) =  (\sqrt{1-\xi_\e^2 (s)} \cos \theta,
\sqrt{1-\xi_\e^2 (s)} \sin \theta, \xi_\e (s));
\label{eq:6.3}
\end{equation}
this is a perturbation of the unit normal $n$, and in fact
\begin{equation}
\bar{n}_\e (s, \theta) - n(s, \theta) = (\chi_\e (s) \,\cos \theta,
\chi_\e (s) \,  \sin \theta, \bar{\chi}_\e (s)).
\label{eq:6.4}
\end{equation}
where $\chi_\e$ and $\bar{\chi}_\e$ are supported in $(-\infty,-s_\e +
2] \cup [s_\e -2,+\infty)$ and satisfy
\begin{equation}
\cosh s\,|\nabla^k \chi_\e| + \cosh^2 s \, |\nabla^k \bar{\chi}_\e |
\leq c_k ,
\label{eq:6.5}
\end{equation}
for all $k \geq 0$.

We now look for all CMC surfaces near the rescaled catenoid which admit
the parametrization
\begin{equation}
{\bf x}_w:[-s_\e, s_\e] \times S^1 \ni (s, \theta) \longrightarrow {\bf
x}^c_\e (s,\theta)+w(s,\theta)\, \bar{n}_\e (s,\theta),
\label{eq:6.6}
\end{equation}
for some smooth, sufficiently small function $w$. By construction, 
these surfaces are normal graphs over $\Sigma^c_\e$ when $|s| \leq s_\e 
-2$ and are vertical graphs when $|s| \geq s_\e -1$. 
The reason for parametrizing surfaces using $\bar{n}_\e$ is so that 
their boundaries are vertical graphs over a circle.  It
follows from the analysis of Appendix II and from (\ref{eq:5.44}) that
such a surface  is CMC if and only if $w$ satisfies a certain nonlinear
equation of the form
\begin{equation}
\ds \frac{1}{\e^2 \cosh^2 s}\,{\cal L}w = \ds \frac{1}{\e^2 \cosh^2 s}\,
\left(-H_0 \, \e^2\, \cosh^2 s + \bar{Q}_\e(w)\right),
\label{eq:6.7}
\end{equation}
where
\[
\begin{array}{rlll}
\bar{Q}_\e(w) & = & \ds L_\e w + \e\, \bar{Q}_\e'\left(\frac{w}{\e \cosh
s}, \frac{\nabla w}{\e \cosh s},\frac{\nabla^2  w}{\e \cosh
s}\right)\\[3mm] & + & \ds \e \cosh s \,\bar{Q}_\e'' \left(\frac{w}{\e
\cosh s}, \frac{\nabla w}{\e \cosh s},\frac{\nabla^2 w}{\e \cosh
s}\right).
\end{array}
\]
Here $\bar{Q}'_\e$ and $\bar{Q}''_\e$ have the same properties
(\ref{eq:5.55}) and (\ref{eq:5.66}) as the functions $Q'_\e$ and
$Q''_\e$ in (\ref{eq:5.44}) and the linear operator $L_\e$ is supported
in $|s| \geq s_\e - 2$, has coefficients of the order  $1/\cosh^{2} s$,
and represents the difference between the Jacobi operator for surfaces
parametrized normally to $\Sigma_\e^c$ and those parametrized using
${\bar{n}}_\e$.   (To see this estimate on the size of the coefficients
of $L_\e$, note from (\ref{eq:1.11}) in Appendix I that this difference
involves $1 -  n \cdot \bar{n}_\e$, which  by (\ref{eq:6.4}) and
(\ref{eq:6.5}) is of order $1/\cosh^2 s$.)
 
Now, given $\phi''_{\pm} \in \pi''\left({\cal
C}^{2,\alpha}(S^1)\right)$,  we wish to solve the boundary value problem
\begin{equation}
\left\{ \begin{array}{rlll} {\cal L} w & = & - \ds H_0 \, \e^2 \, \cosh^2 s +
\bar{Q}_\e ( w ) \qquad &  \mbox{in} \qquad ( -s_\e, s_\e ) \times S^1
\\[3mm] \pi'' w  & = & \phi''_{\pm} \qquad & \mbox{on} \qquad \{\pm s_\e
\}  \times S^1 .
\end{array}
\right.
\label{eq:6.8}
\end{equation}
A solution will produce a CMC surface with boundary components
parametrized by
\[
S^1 \ni \theta \longrightarrow \left(\e \cosh s_\e \, \cos \theta, \e
\cosh s_\e \, \sin \theta, \pm \e s_\e - \e \log \e + \e \log 2  \mp
w(\pm s_\e, \theta) \right).
\]
Note these are vertical graphs over (small) circles.

We solve (\ref{eq:6.8}) by a standard contraction mapping argument.
First fix $\delta \in (1,2)$ and define
\begin{equation}
\tilde{w} = {\cal P}_{s_\e}(\phi''_\pm)- H_0 \, {\cal G}_{s_\e}(\e^2\, \cosh^2 s)
\label{eq:6.9}
\end{equation}
as an approximate solution for the problem. Then, writing  $w =
\tilde{w} + v$, we must find a function $v \in  {\mathcal
C}^{2,\alpha}_{\delta}([-s_\e, s_\e]  \times S^1)$ such that
\begin{equation}
\left\{ \begin{array}{rlll} {\cal L}v & = & \ds \bar{Q}_\e(\tilde{w}+v)
\qquad & \mbox{in} \qquad (-s_\e,s_\e)\times S^1 \\[3mm] \pi''v & = & 0
\qquad & \mbox{on} \qquad \{\pm s_\e \} \times S^1.
\end{array}
\right.
\label{eq:6.12}
\end{equation}
This will be accomplished by finding a fixed point of the mapping
\begin{equation}
{\cal N}_\e (v) \equiv {\cal G}_{s_\e }\left(\bar{Q}_\e (\tilde{w} + v)
\right).
\label{eq:6.13}
\end{equation}
Although not explicit in the notation, this operator depends on
$\phi''_\pm$.
\begin{proposition} 
Fix $\kappa >0$. Then there exist constants $c_\kappa >0$ and $\e_0 >0$
such that if $0 < \e \leq \e_0$ and if $\phi''_{\pm} \in
\pi''\left({\cal C}^{2,\alpha} (S^1)\right)$ is fixed with
$||\phi''_{\pm}||_{2,\alpha} \leq \kappa \, \e^{3/2}$,  then ${\cal
N}_\e$ is a contraction mapping on the ball
\[
B_{c_\kappa} \equiv \{v: ||v||_{2,\alpha,\delta} \leq  2\, c_\kappa
\,\e^{(8 +\delta)/4} \},
\]
and hence has a unique fixed point in this ball.
\label{pr:9}
\end{proposition}
{\bf Proof:} We must show that
\[
||{\cal N}_\e (0)||_{2,\alpha,\delta}\leq c_\kappa\, \e^{(8+\delta)/4}
\]
and
\[
||{\cal N}_\e (v_2)-{\cal N}_\e (v_1)|| _{2,\alpha, \delta} \leq
\frac{1}{2} ||v_2-v_1||_{2,\alpha, \delta },
\]
for all $v_1, \, v_2 \in B_{c_\kappa}$. For then, if $v \in B_{c_\kappa}$,
then 
\[
||N(v)||_{2,\alpha,\delta} \leq ||N(0)||_{2,\alpha,\delta} +
||N(v) - N(0)||_{2,\alpha,\delta} \leq 2c_\kappa \e^{(8+\delta)/4}.
\]

We begin with the first of these. To do this, we must estimate
$\tilde{w}$.  Set $\tilde{w}_0 =  {\cal P}^0_{s_\e}(\phi''_\pm)$; then
since $e^{-2s_\e} = \e^{1/2}$ and $||\phi''_{\pm}||_{2,\alpha} \leq
\kappa  \,\e^{3/2}$, we get from Proposition~\ref{pr:88} and
Proposition~\ref{pr:8} that
\begin{equation}
||\tilde{w}_0||_{2,\alpha, 2} + ||{\cal P}_{s_\e}(\phi''_\pm)
-\tilde{w}_0||_{2,\alpha,\delta} \leq c \, \kappa \, \e^{2}.
\label{eq:6.10}
\end{equation}
Next, even if the result of Proposition~\ref{pr:6} does not hold when
the weight parameter $\delta =2$ and taking advantage of the fact that
$\e^2\, \cosh^2 s$ is independent of $\theta$, we can use directly
(\ref{eq:5.6}), to estimate
\begin{equation}
||{\cal G}_{s_\e}(\e^2\, \cosh^2 s) ||_{2, \alpha, 2} \leq c \,\e^2.
\label{eq:6.11}
\end{equation}
Notice that $\|\cdot \|_{2, \alpha, 2 }\leq \| \cdot \|_{2, \alpha,
\delta}$, since $\delta \in (1,2)$.  Putting these
together, we get
\begin{equation}
||\tilde{w}||_{2,\alpha,2} \leq c \,  \e^2,
\label{eq:6.111}
\end{equation}
for some constant $c$ depending on $\kappa$ but independent of $\e$.

Next we estimate $||{\cal N}_\e (0)||_{2,\alpha, \delta}$ by
\[
||{\cal N}_\e (0)||_{2,\alpha, \delta} \leq c \, (  \|L_\e
\tilde{w}\|_{0,\alpha,\delta} + \|\e\, \bar{Q}'_\e (\tilde{w}/\e \cosh
s) \|_{0,\alpha, \delta} + \|\e \cosh s \, \bar{Q}''_\e (\tilde{w}/\e
\cosh s) \|_{0,\alpha,\delta}) .
\]
We have first
\[
\left\|L_\e \tilde{w}\right\|_{0,\alpha,\delta} \leq c \, \e^{(8 +
\delta)/4},
\]
and then
\[
\left\|\e\, \bar{Q}'_\e \left(\frac{\tilde{w}}{\e \cosh s} \right)
\right\|_{0,\alpha, \delta} \leq  c \, \e^{(10 + \delta)/4} \qquad
\mbox{\rm and} \quad \left\| \e \, \cosh s \, \bar{Q}''_\e \left(
\frac{\tilde{w}}{\e \cosh s}  \right) \right\|_{0,\alpha, \delta} \leq
c \,\e^{(12+ \delta)/4}.
\]
Again, all constants depend on $\kappa$ but not on $\e$.

Now clearly it suffices to choose $\e$ sufficiently small and $c_\kappa$
equal to twice the constant in (\ref{eq:6.111}) in order  for the stated
estimate for ${\cal N}_\e (0)$ to hold.

For the other estimate, if $v_1,v_2 \in B_{c_\kappa}$, then
\[
\left\| {\cal G}_{s_\e}L_\e(v_1 - v_2)\right\|_{2,\alpha,\delta} \leq c
\, \left\| L_\e (v_2 -v_1) \right\|_{0,\alpha,\delta}  \leq c \,
\e^{1/2} \, ||v_2 -v_1||_{2,\alpha,\delta},
\]
\[
\left\| \e  \,\left( \bar{Q}'_\e \left( \frac{\tilde{w}+v_2}{\e \cosh s}
\right)- \bar{Q}'_\e \left(\frac{\tilde{w}+v_1}{\e \cosh s } \right)
\right) \right\|_{0,\alpha, \delta} \leq  c \, \e \,  ||v_2-v_1||_{2,
\alpha, \delta}
\]
and finally
\[
\left\| \e \cosh s \,\left( \bar{Q}''_\e \left( \frac{\tilde{w}+v_2}{\e
\cosh s} \right)  -   \bar{Q}''_\e \left(\frac{\tilde{w}+v_1}{\e \cosh s
} \right)  \right)  \right\|_{0,\alpha, \delta} \leq  c\, \e^{3/2}\,
||v_2-v_1||_{2, \alpha, \delta}.
\]
We are using here that $||\tilde{w}/\cosh s||_{2,\alpha,0} \leq c \,
\e^{7/4}$ and $||v_i/\cosh s||_{2,\alpha,0} \leq c \, \e^{9/4}$. The
result follows at once for all $\e$ small enough. \hfill $\Box$

\ms

To conclude this section we examine the Cauchy data map
\begin{equation}
{\cal S}_\e : \left( \pi'' \left( {\cal C}^{2,\alpha}(S^1) \right)
 \right)^2  \longrightarrow  \left( {\cal C}^{2,\alpha}  (S^1)\times
 {\cal C}^{1,\alpha} (S^1)\right)^2
\label{eq:cdw}
\end{equation}
given by
\[
\begin{array}{rllll}
{\cal S}_\e(\phi''_{\pm}) & = & \left((\e \, s_\e + \e \log (2/\e) - w (s_\e, \cdot), \e  -  \partial_s
w (s_\e, \cdot))\right. , \\[3mm]
&   &  \, \, \, \left. (-\e \, s_\e + \e \log (2/\e) + w (- s_\e, \cdot) , -
\e  - \partial_s w (- s_\e, \cdot) \right)
\end{array}
\]
where $w$ is the solution of (\ref{eq:6.12}) given by
Proposition~\ref{pr:9}.  We shall also consider the Cauchy data map
${\cal S}_0$ for the operator ${\cal P}_{s_\e}^0$. It is simple to check
that
\[
{\cal S}_0(\phi''_{\pm}) = \left( (\e \, s_\e + \e \log(2/\e) -
\phi''_{+}, \e - |D_\theta| \phi_{+}'' ), (-\e \, s_\e +\e \log (2/\e) +
\phi''_{-} , -\e + |D_\theta| \phi_{-}'' )\right).
\]

The comparison between these two Cauchy data mappings plays a key role
in our construction.

\begin{corollary}
The mappings ${\cal S}_\e$ and ${\cal S}_0$ are continuous. Furthermore,
there exists a constant $c > 0$ such that for any $\kappa > 0$ there
exists an $\e_0 > 0$ such that if $\e \in (0, \e_0]$, then for all
$\|\phi''_{\pm}\|_{2,\alpha} \leq \kappa \, \e^{3/2}$, we have
\begin{equation}
||({\cal S}_\e -{\cal S}_0)(\phi''_{\pm})||_{({\cal C}^{2,\alpha}
\times {\cal C}^{1, \alpha})^2 } \leq c \, \e^{3/2}.
\label{eq:6.14}
\end{equation}
\label{co:1}
\end{corollary}
{\bf Proof}: The statement about continuity is straightforward and is
left  to the reader.  Next, we must estimate the Cauchy data for the
function $w - {\cal P}_{s_\e}^0(\phi''_\pm)$. By Proposition~\ref{pr:8},
the Cauchy data of the function ${\cal P}_{s_\e}^0(\phi''_\pm)$ differs
from that of the function $w_\e = {\cal P}_{s_\e}(\phi''_\pm)$ by a term
of order $e^{-2s_\e}||\phi''_\pm||_{2,\alpha} \leq c \e^2$.  Therefore
we must estimate the Cauchy data of the function $w - w_\e = \tilde{w} +
v - w_\e$. Now it suffices to use  (\ref{eq:6.111}) and the fact that
$||v||_{2,\alpha,\delta} \leq  2\,c_\kappa\, \e^{\delta/4} \, \e^{2}$.
This ends the proof of the Corollary. \hfill $\Box$

\ms

It is important here that the constant $c$ is independent of $\kappa$.
(See Definition 3 in \S 8 for the precise meaning with which this is
meant to be understood.)

\section{The mean curvature operator for a graph}
\label{sse:3.1}

Assume that $\Sigma \subset {\R}^3$ is a regular surface such that $0
\in \Sigma$  and $T_0\Sigma$ is the $x\,y$-plane. Then $\Sigma$ can be
locally  parameterized, near the origin, as a graph
\begin{equation}
{\bx} :  B_{\bar{\rho}} \ni (x,y) \longrightarrow (x,y, u (x,y)) \in
\Sigma \subset {\R}^3,
\label{eq:1.1}
\end{equation}
where $u : B_{\bar{\rho}}  \longrightarrow {\R}$ is a regular function
which  satisfies
\begin{equation}
\nabla^{k} u (x,y)  = {\cal O} ( r^{2-k}), \quad k \leq 2,  \qquad
\nabla^{k} u (x,y)  = {\cal O} (1), \quad k \geq 3,
\label{eq:1.2}
\end{equation}
where $r= (x^2+y^2)^{1/2}$. In this parameterization, the mean
curvature operator $H_{u} (x,y)$ of the graph defined by the function
$u$ at  the point of parameter $(x,y)$ is then given by \cite{GT}
\begin{equation}
H_{u} (x,y) = \displaystyle \nabla \left( \frac{\nabla u}{(1 +|\nabla
u|^2)^{1/2}} \right).
\label{eq:1.3}
\end{equation}
Notice that we have defined the mean curvature to be the sum of the
principal  curvatures $H = k_1 +k_2 $, not the average.

By our assumptions, all surfaces sufficiently close to $\Sigma$ can  be
parameterized, in some neighborhood of $0$, as a vertical graph over  a
neighbourhood of $0$ in the $x\,y$-plane, namely as
\begin{equation}
B_\rho \ni (x,y) \longrightarrow (x, y, u (x,y) + w (x,y))
\label{eq:1.4}
\end{equation}
for some (regular) function $w : B_\rho \longrightarrow {\R}$.

It follows from (\ref{eq:1.3}) that the linearized mean curvature
operator about $\Sigma$ is given explicitly by
\begin{equation}
\Lambda_u : w \longrightarrow  \displaystyle  \nabla \left( \frac{\nabla
w}{ (1+|\nabla u |^2)^{1/2}} - \frac{\nabla u \cdot \nabla w}{(1+|\nabla
u |^2)^{3/2}} \nabla u  \right).
\label{eq:1.5}
\end{equation}
Performing the change of variable $(x,y) = e^{-t} (\cos \theta, \sin
\theta)$,  the linearized operator, which we still denote by
$\Lambda_u$, has the form
\begin{equation}
\Lambda_u : w \longrightarrow  \displaystyle e^{2t}  \,(\del^2_{t t}
+\del^2_{\theta \theta}) + {\Lambda}'_u,
\label{eq:1.6}
\end{equation}
where ${\Lambda}'_u $ is a second order partial differential operator
with no terms of order zero and with coefficients bounded in ${\cal
C}^\infty  ([-\log \rho , +\infty) \times S^1)$.

We also may expand the mean curvature $H_{u+w}$ of the surface
parameterized  by (\ref{eq:1.4}) in terms of the mean curvature $H_u$ of
$\Sigma$ and  $\Lambda_u$; thus
\begin{equation}
H_{u+w}  = H_{u} + \Lambda_u w -  Q'_u (e^t \, \nabla w , e^t \,
\nabla^2 w)  - e^t \, Q''_u   (e^t \, \nabla w , e^t \, \nabla^2 w),
\label{eq:1.8}
\end{equation}
where $Q'_u$ and $Q''_u$ are functions with coefficients bounded in
${\cal C}^\infty ([-\log \rho, +\infty) \times S^1)$ which satisfy
\[
Q'_u (0,0)= Q''_u (0,0) = 0 \qquad \mbox{and}\qquad \nabla Q'_u (0,0)=
\nabla Q''_u (0,0) =0 .
\]
and also
\[
\nabla^2 Q''_u (0,0)=0.
\]
These facts are established in Appendix II.

\section{Analytic modification of a surface $\Sigma_0 \subset {\R}^3$ using 
Green's function}

\label{se:4}

Assume that $\Sigma_0$ is a regular, orientable CMC surface with boundary,
positioned and oriented as in the previous section. We use a local
chart $\bx$ as in (\ref{eq:1.1}), with $\bar{\rho} < 1$.   We also
assume that (\ref{eq:1.2}) holds in $B_{\bar{\rho}}$.  We define, for
$\rho \leq \bar{\rho}/2$,
\[
\Sigma_0 (\rho) \equiv \Sigma_0 \setminus {\bx} \left( B_\rho \right).
\]
Still assuming that $\rho \leq \bar{\rho}/2$, in $\Sigma_0$, we choose a
unit  vector field $\bar{\nu}$ which is equal to a normal unit vector
field $\nu$  in $\Sigma_0 (2 \rho)$ and which is equal to $(0,0,1)$ in
${\bx}(B_\rho)$.  We assume that $\nu \cdot \bar{\nu}  \geq 1/2$ on all
$\Sigma_0$.  All surfaces near to $\Sigma_0$ can be parameterized by
$\Sigma_0 \ni p \longrightarrow p + w(p)\, \bar{\nu}(p)$ for some  small
function $w$. The linearized mean curvature operator
\[
\Lambda: {\cal C}_{\cal D}^{2, \alpha}(\Sigma_0)= \{ w \in {\cal C}^{2,
\alpha} (\Sigma_0): w=0 \quad \mbox{on} \quad \partial \Sigma_0\}
\longmapsto  {\cal C}^{0, \alpha}(\Sigma_0),
\]
relative to this vector field has the familiar form
\[
\Lambda \equiv \Delta_{\Sigma_0} + |{\bf A}_{\Sigma_0}|^2
\]
in $\Sigma_0(2\rho)$, while in ${\bx} (B_\rho)$ it is given by 
(\ref{eq:1.5}).  Although not obvious at this stage, the use of 
$\bar{\nu}$ to parametrize nearby surfaces is intended to make
the later analysis simpler.

By construction, $\Lambda$ depends on $\rho$. It follows from the
analysis in Appendix I, cf. particularly (\ref{eq:1.11}), that 
$\Lambda$ tends to $\Delta_{\Sigma_0} + |{\bf A}_{\Sigma_0}|^2$ as
$\rho \to 0$. In particular, if $\Sigma_0$ is nondegenerate in the
sense of Definition \ref{de:1}, then $\Lambda$ is an isomorphism for 
$\rho$ sufficiently small. From now on, we shall  assume that $\rho \leq
\bar{\rho}/2$ is fixed once for all so that this is the case. We 
may then solve the equation
\begin{equation}
\Lambda \gamma_0 = - 2 \pi \, \delta_0, \qquad \mbox{in} \qquad \Sigma_0,
\label{eq:2.2}
\end{equation}
with $\gamma_0 = 0$ on $\partial \Sigma_0$.

The following Lemma follows easily from (\ref{eq:1.5}), 
using (\ref{eq:1.2}), and details are left to the reader.
\begin{lemma}
Assume that (\ref{eq:1.2}) holds and that $\gamma_0$ is the solution of
(\ref{eq:2.2}). Then there exist constants $a_0, a_1, a_2 \in {\R}$
such that, for all $k \geq 0$,
\begin{equation}
\nabla^k \left( \gamma_0 (x,y) - \bar{\gamma}_0 (x,y)\right)  =  {\cal
O} (r^{2-k} \, \log 1/r),
\label{eq:2.3}
\end{equation}
where $\bar{\gamma}_0 (x,y) = -\log r + a_0 + a_1 \, x + a_2 \, y $.
\label{le:1}
\end{lemma}

For $0 < \e$ we define the surface $\bar{\Sigma}_\e$ to be the 
one parameterized by 
\begin{equation}
\Sigma_0 \setminus B_\e \ni p \longrightarrow p + \e \, \gamma_0 (p)
\, \bar{\nu}(p) \in {\R}^3.
\label{eq:2.4}
\end{equation}
If $\e$ is small enough, this is a regular surface. 

We now compare the mean curvatures of $\bar{\Sigma}_\e$ and $\Sigma_0$.
\begin{proposition}
We may estimate the difference between $H_\e$, the mean curvature of 
$\bar{\Sigma}_\e$, and $H_0$, the mean curvature of $\Sigma_0$, by
\[
\nabla^k \left(H_\e - H_0\right) = {\cal O}\left(r^{-k}(\e^2\,
r^{-2}+\e^3\, r^{-4})\right)\qquad \mbox{in} \qquad \bx (B_\rho
\setminus B_\e)
\]
and by
\[
\nabla^k  \left(  H_\e  - H_0 \right) =  {\cal O} \left(\e^2 \right)
\qquad \mbox{in} \qquad \Sigma_0 (\rho)
\]
for all $k \geq 0$. 
\label{pr:4}
\end{proposition}
{\bf Proof :}  This follows at once from (\ref{eq:1.8}) with $w=
\e \, \gamma_0$. \hfill $\Box$

\begin{corollary} $H_\e$ is bounded independently of $\e$ in 
$\Sigma_0 \setminus B_{c\e^{3/4}}$ for any fixed constant $c>0$.
\label{co:mcbds}
\end{corollary}

Now, from Lemma~\ref{le:1}, if $\bar{\Sigma}_\e$ is translated vertically 
(along the $z$-axis) by $- \e \, a_0$, then it will parameterized 
near $0$ by
\begin{equation}
\bx_\e: \Sigma_0 \setminus B_\e \ni (x,y) \longrightarrow 
\left( x, y , - \e \, \log r + u (x,y) +  {\cal O} (\e \, r) \right).
\label{eq:4.1}
\end{equation}
Comparing this expansion with the one in (\ref{eq:3.4}) and using that $u
(x,y) = {\cal O} (r^2)$, we see that the vertical distance between the two
surfaces is estimated by ${\cal O}\left(r^2 +\e\, r+\e^3\, r^{-2}\right)$.  
We have chosen the vertical translations of both the catenoid and $\Sigma_0$ 
carefully to minimize the distance between them. At any rate, this 
quantity is minimal for $r \sim \e^{3/4}$.  This and 
Corollary~\ref{co:mcbds} make it now quite reasonable that 
we restrict attention to a neighbourhood $r \geq c \e^{3/4}$, for any
$c>0$. Thus we now define $\Sigma_\e$ to be the surface which
is given near the origin by the parametrization
\[
\Sigma_0 (c \, \e^{3/4}) \ni p \longrightarrow p + \e \,  \gamma_0 (p)
\, \bar{\nu}(p) \in {\R}^3.
\]
The constant $0 < c < 1/8$ is now fixed once and for all.
Notice that on $\Sigma_\e$, $\e\, r={\cal O}(r^2)$. The `inner
boundary' of $\Sigma_\e$, created by the excised ball, will
be denoted $\del_1\Sigma_\e$. 

\section{Geometric modifications of the surface $\Sigma_\e$}
\label{se:6}

In order to match the Cauchy data of perturbations of the catenoid
and of $\Sigma_\e$, it is necessary to allow some extra
flexibility in the boundary data, specifically in the low
($j = 0, \pm 1$) eigenmodes. This flexibility is created by
considering not just the surface $\Sigma_\e$, but also a family
of modifications of it, comprised of rotations and translations 
and alterations of the parameter $\e$; the effect of these
modifications in the boundary data is seen only in the low 
eigenmodes. In this section we define and analyze this family.

The parameter set for this family of surfaces will be
denoted ${\cal A}= ((T_1, T_2,  T_3), (R_1,R_2), e) \in 
\cal U$, where ${\cal U}$ is a neighbourhood of the
origin in ${\R}^{3} \times {\R}^{2} \times {\R}$.
The effect of the parameters $(T_1,T_2,T_3)$ will be to
translate the surface by this vector. The parameters $(R_1,R_2)$
correspond to a rotation of the surface by the matrix 
\[
\exp \left( \begin{array}{ccc} 0   & 0 & - R_1 \\[3mm] 
0  & 0  & - R_2 \\[3mm] R_1  & R_2  & 0 \end{array}\right ),
\]
which has the form
\[
\left( \begin{array}{ccc} 1 + {\cal O}(|R|^2)   & {\cal O}(|R|^2)
& - R_1 +{\cal O}(|R|^3)\\[3mm] {\cal O}(|R|^2)       & 1 + {\cal
O}(|R|^2)   & - R_2 +{\cal O}(|R|^3)\\[3mm] R_1 + {\cal O}(|R|^3) & R_2
+ {\cal O}(|R|^3) & 1 + {\cal O}(|R|^2)
\end{array}
\right ).
\]
Finally, $e$ changes the scaling parameter $\e$ into $\e - e$.
Since these operations do not commute, we make the convention
that ${\cal A}$ acts on $\Sigma_\e$ by first changing $\e$ to
$\e - e$, next translating by $T_3$ in the vertical direction, 
then performing the rotation and finally translating by
$(T_1,T_2)$ horizontally.

If the neighbourhood ${\cal U}$ is sufficiently small, the 
resulting surface, which will be denoted ${\Sigma}_{\e, {\cal A}}$, 
can still be locally parameterized as a graph over the $x \,y$-plane. 
We shall define a norm for the vector ${\cal A}$ by
\begin{equation}
\| {\cal A} \| \equiv \e^{1/4} \,  \| (T_1, T_2)\|_{\R^2} + |\log \e
|^{-1} \,  |T_3| + \e^{3/4} \, \| (R_1, R_2)\|_{\R^2} + |e|.
\label{eq:4.2}
\end{equation}
This choice of scaling factors on the various components of
${\cal A}$ is necessitated by the analytic details of the
ensuing arguments. In fact, this norm is related to the function $S^1 
\ni \theta \rightarrow (\xi (\theta), r \, \del_r\xi (\theta))$ for 
$r\sim \e^{3/4}$, where
\[
\xi : S^1 \ni \theta \longrightarrow e \, \log r + T_3 + r \, (R_1
\cos \theta + R_2 \sin\theta) + \e \, r^{-1}(T_1 \cos \theta + T_2
\sin \theta).
\]
Hence it measures the effect of the geometric modifications 
on the set of points where the gluing will be done, see
Proposition 10.

As we have already noted, in some neighbourhood of its inner boundary, 
$\Sigma_{\e,{\cal A}}$ can be written as a graph over the $x\, y$-plane,
and this graph function can be compared to the graph function
for the original surface $\Sigma_0$ and also to the one for the 
catenoidal (or rather, logarithmic) end. These `comparison' graph 
functions will be 
denoted $w_m$ and $\hat{w}_m$, respectively. These functions depend 
on all the parameters. The main result of this section gives estimates 
on these functions, but first we introduce some convenient
notation.
\begin{definition}
Henceforth the notation $f = {\cal O}(g(\e,r))$ shall mean that
the function $f$ (usually on $\Sigma_{\e,{\cal A}}$) is bounded
by a constant $c$ times the function $g$ of $\e$ and $r$, i.e.
$f \leq c\, g(\e,r)$, where the constant $c$ does not depend
on either $\e$ or $\kappa$. On the other hand, $f = {\cal O}_\kappa
(g(\e,r))$ shall mean that $f$ is bounded similarly, but
by a constant $c_\kappa$ which is allowed to depend on
$\kappa$, but is still independent of $\e$. Furthermore,
a bound of the latter type may sometimes by converted to
a bound of the former type as follows. If, for example,
$f = {\cal O}_\kappa(\e^{3/2})$, and if we have (as shall
always be true in all the calculations below) that
the constant $c_\kappa$ in this estimate is bounded by
a fixed polynomial in $\kappa$, then we also have $f = 
{\cal O}(\e^{3/4})$, since $c_\kappa \e^{3/2} \leq
c \, \e^{3/4}$ provided $\e$ is sufficiently small, for
fixed $\kappa$. This reasoning will be justified because,
although we need the flexibility to set $\kappa$ fairly
large, once we have done so it will be fixed, and this will
then determine an upper bound $\e_0$ for $\e$.
\label{de:constants}
\end{definition} 

\begin{proposition}
Fix $\kappa >0$. Then there exists an $\e_0 >0$ such that if 
$0 < \e \leq \e_0$ and $\|{\cal A}\|\leq \kappa \, \e^{3/2}$, 
then for $\frac{1}{4}\e^{3/4} \leq r \leq \rho$, the surface 
${\Sigma}_{\e,{\cal A}}$ can be parameterized as
\begin{enumerate}
\item 
\begin{equation}
(x,y) \longrightarrow \left(x,y, u (x,y) + w_m (x,y) \right),
\label{eq:4.3}
\end{equation}
where $w_m(x,y)={\cal O}_\kappa \left (\e^{3/4}\, r + \e\,|\log r|
\right)$ and in addition, for all $k \geq 1$, $\nabla^k
w_m(x,y)={\cal O}_\kappa \left(r^{-k}\left(\e^{3/4}\, r +\e
\right)\right)$. 
\item 
\begin{equation}
(x,y)\longrightarrow \left(x,y,-\e\, \log r +\hat{w}_m(x,y) \right),
\label{eq:4.4}
\end{equation}
where $\hat{w}_m (x,y)={\cal O}_\kappa \left(r^{2}+\e^{3/2}\,|\log \e| 
\right)$ and for $k\geq 1$, $\nabla^k \hat{w}_m (x,y) = 
{\cal O}_\kappa \left(r^{2-k} + \e^{3/2} r^{-k}\right)$.
\end{enumerate}
\label{pr:5}
\end{proposition}

\noindent
{\bf Proof :} We shall only prove the estimates for ${w}_m$ and 
$\hat{w}_m$ because the estimates for the derivatives follow
from these in a straightforward manner.

First, notice that from $\|{\cal A}\| \leq \kappa \, \e^{3/2}$ we have 
\[
\|(T_1, T_2)\|_{\R^2} \leq \kappa\, \e^{5/4}, \qquad \qquad |T_3| \leq
\kappa \, \e^{3/2} \,|\log \e | ,
\]
\[
|e|\leq \kappa \, \e^{3/2}  , \qquad \qquad \| (R_1, R_2)\|_{\R^2}
\leq \kappa \, \e^{3/4}.
\]

We perform the transformation on all of ${\R}^3$, first translating
vertically by $-e \, (\gamma_0 - a_0) + T_3 = e\log r + T_3 +
{\cal O}(er)$, then applying the rotation matrix 
and finally translating horizontally by the vector $T' = (T_1,T_2)$. 
Acting on all of space, this effects a change from the coordinates
$(x,y,z)$ to the coordinates $(\tilde{x}, \tilde{y}, \tilde{z})$. The
precise relationship is
\[
\begin{array}{rlll}
\tilde{x}  & = & \displaystyle T_1 + (1 + {\cal O}(|R|^2)) x + {\cal O}
(|R|^2) y  - \displaystyle (R_1+{\cal O}(|R|^3)) \left (e\,\log r + T_3
+ z + {\cal O}(e r)\right),\\[3mm] 
\tilde{y} & = & T_2 + (1 + {\cal O}(|R|^2)) y +{\cal O}(|R|^2) x - \displaystyle (R_2+{\cal O}(|R|^3))
\left( e \, \log r + T_3 + z + {\cal O}(e r)\right), \\[3mm]
\tilde{z} & = & (R_1 + {\cal O}(|R|^2)) x +(R_2 + {\cal O}(|R|^2) y 
+ (1 + {\cal O}(|R|^2) )\left(e \, \log r + T_3 + z + {\cal O}(e r) \right).
\end{array}
\]
Recalling that $z = \e \, (\gamma_0(x,y) - a_0) + u(x,y)$, we first observe that
\[
|(\tilde{x},\tilde{y}) - (x,y)| = {\cal O}
(|R|^2 r + |R| ( \e \,|\log r| + |T_3| + r^2 + \e r) + |T'|) = 
{\cal O}_\kappa(\e^{3/4}r^2 + \e^{5/4}),
\]
for $\e$ sufficiently small.  Hence if we set $\tilde{r} =
|(\tilde{x},\tilde{y})|$, we obtain
\[
\tilde{r} = r + {\cal O}_\kappa(\e^{3/4} r^2 + \e^{5/4}) , 
\]
and in particular, we get for all $r \geq \e^{3/4}/8$, we can state that 
$r/2 \leq \tilde{r}\leq 2 r $ provided $\e$ is small enough. Using this 
first information, we obtain
\begin{equation}
\log r = \log \tilde{r} + {\cal O}_\kappa(\e^{3/4} \tilde{r} + \e^{5/4}\tilde{r}^{-1})
\label{eq:logrrt}
\end{equation}
and also
\begin{equation}
r = \tilde{r} + {\cal O}_\kappa(\e^{3/4} \tilde{r}^2 + \e^{5/4}) .
\label{eq:rrt}
\end{equation}
Inserting these estimates into the equations for $\tilde{x}$
and $\tilde{y}$ above, we see that
\[
(\tilde{x} - x, \tilde{y} - y) = T' + {\cal O}_\kappa 
\left(\e^{7/4}\, |\log \tilde{r} | + \e^{3/4}\, \tilde{r}^2 
\right) = {\cal O}_\kappa \left(\e^{3/4}\, \tilde{r}^2 
+ \e^{5/4}\right).
\]
Thus, we can evaluate
\[
u(x,y) = u(\tilde{x}, \tilde{y}) + \nabla u(\xi,\eta)(x-\tilde{x},y-\tilde{y}) + 
{\cal O}_\kappa (\e^{5/2}+ \e^{3/2}\tilde{r}^4),
\]
where $(\xi,\eta)$ is some point on the line between $(x,y)$ and
$(\tilde{x},\tilde{y})$. Since $|(\xi,\eta)| = {\cal O}(r)$, 
we obtain $|\nabla u(\xi,\eta)| = {\cal O}(r)$, hence
\begin{equation}
u(x,y) = u(\tilde{x}, \tilde{y}) + {\cal O}_\kappa (\e^{3/4} \tilde{r}^3+ \e^{5/4}\tilde{r}).
\label{eq:ldt}
\end{equation}

It is also an easy matter to check that
\[
 \tilde{z} = z + {\cal O}_\kappa(\e^{3/4}\tilde{r} + \e^{3/2}|\log \e|).
\]
and also that
\[
u(x,y) = u(\tilde{x}, \tilde{y}) + \nabla u(\xi,\eta)(x-\tilde{x},y-\tilde{y}) + 
{\cal O}_\kappa (\e^{5/2}+ \e^{3/2}\tilde{r}^4) ,
\]
where $(\xi,\eta)$ is some point on the line between $(x,y)$ and
$(\tilde{x},\tilde{y})$. Since $|(\xi,\eta)| = {\cal O}(r)$, 
we obtain $|\nabla u(\xi,\eta)| = {\cal O}(r)$.
Finally, recalling again that $z = \e (\gamma_0(x,y) - a_0) + u(x,y)$, and
collecting the previous estimates, we get
\[
\tilde{z} = u(\tilde{x},\tilde{y}) + {\cal O}_\kappa(\e^{3/4}\tilde{r} + \e |\log \tilde{r}|),
\]
which gives the desired estimate for $w_m$. 

For the other part of the proposition, we wish to estimate the
function $\hat{w}_m$, where
\[
\tilde{z} = -\e \log \tilde{r} + \hat{w}_m(\tilde{x},\tilde{y}).
\]
This time, we use the estimate (\ref{eq:logrrt}) to get
\[
\hat{w}_m = {\cal O}_\kappa(\tilde{r}^2 + \e^{3/2} |\log \e|).
\]
\hfill $\Box$

\begin{corollary}
The mean curvature $H_{\e,{\cal A}}$ of $\Sigma_{\e,{\cal A}}$
satisfies the same estimates as that of $\Sigma_\e$, namely
\[
|\nabla^j (H_{\e,{\cal A}} - H_0)| = {\cal O} (r^{-j}(\e^2 
r^{-2} + \e^3 r^{-4})).
\]
\label{co:estmcea}
\end{corollary}
{\bf Proof:} Because $\Sigma_{\e,{\cal A}}$ is obtained from
$\Sigma_\e$ by first modifying the Green function by
an amount much less than $\e$ and then applying a rigid
motion, it is clear that $H_{\e,{\cal A}} - H_0$ and all
its derivatives are bounded by a multiple of $\e^2$ 
outside $B_{\bar{\rho}}$. Inside this ball we know that
$|H_{\e,{\cal A}} - H_0| = {\cal O} (\e^2 r^{-2}
+ \e^3 r^{-3})$.  From (\ref{eq:rrt}) one easily obtains
\[
\e^2 r^{-2} + \e^3 r^{-4} = {\cal O} (\e^2 
\tilde{r}^{-2} + \e^3 \tilde{r}^{-4} ),
\]
as desired. The bounds for the derivatives are similar. \hfill $\Box$

\ms
We also require the following result.
\begin{proposition} 
If $r \in [\frac{1}{4}\e^{3/4}, 4\, \e^{3/4}]$, then the 
parameterization of ${\Sigma}_{\e, {\cal A}}$ has an expansion
of the form
\begin{equation}
(x,y) \longrightarrow (x,y, -\e \log r + w_m^0(x,y) + \bar{w}_m(x,y)),
\label{eq:4.5}
\end{equation}
where 
\[
w_m^0(x,y) =  \left(e \log r + T_3 + R_1 x + R_2 y + 
\e \, r^{-2}(T_1 x + T_2 y) \right)
\]
and, for all $k\geq 0$, $|\nabla^k \bar{w}_m(x,y)| = 
{\cal O}(\e^{(6-3k)/4})$. 
\label{pr:5.5}
\end{proposition}
{\bf Proof :} We know, first of all, that 
\[
(\tilde{x} - x, \tilde{y} - y) = T' + {\cal O}_\kappa \left(\e^{7/4}\, |\log \e| \right)
\]
and then that
\[
u(x,y) = {\cal O}(\tilde{r}^2) = {\cal O}(\e^{3/2}).
\]
In both of these we have used the upper bound on $r$. 
Finally, 
\[
r^2 = \tilde{r}^2 - 2(T_1 \tilde{x} + T_2 \tilde{y} ) 
+ {\cal O}_\kappa(\e^{5/2} \, |\log \e|),
\]
so that 
\[
\log r = \log \tilde{r} - \frac{T_1 \tilde{x} +  T_2 \tilde{y}}{\tilde{r}^2} + {\cal O}_\kappa(
\e \, |\log \e|).
\]
Putting these all together in the expression for $\tilde{z}$ yields 
\[
\tilde{z} = -\e \log \tilde{r} + w_m^0(\tilde{x},\tilde{y}) + {\cal O}(\e^{3/2}) +
{\cal O}_\kappa(\e^{7/4}),
\]
which gives the estimate for $\bar{w}_m$. The bounds for its derivatives
are handled similarly.  \hfill $\Box$

\section{The linearized mean curvature operator about 
${\Sigma}_{\e, {\cal A}}$}
\label{se:9}

In this section we shall study the Jacobi operator $\Lambda_{\e,
{\cal A}}$ (relative to a transverse, but not everywhere normal unit 
vector field $\tilde{\nu}$) for the surface $\Sigma_{\e,{\cal A}}$. 
The results we obtain are the usual ones, namely solvability of
$\Lambda_{\e,{\cal A}}u= f$ in appropriate weighted spaces 
with homogeneous Dirichlet boundary conditions (off the low eigenmodes)
as well as of $\Lambda_{\e,{\cal A}}u= 0$ with inhomogeneous
Dirichlet boundary conditions. 

It is slightly simpler to use a different parameterization now. 
Thus let $t = -\log r $ and set  
\begin{equation}
\bar{t} = - \log  \rho + 1 >0 \qquad \mbox{\rm and} \qquad
t_\e  = -\log \left(  \e \cosh (1/4 \,\log \e) \right).
\label{eq:7.1}
\end{equation}
Notice that $t_\e = -3/4 \log \e + \log 2+ {\cal O} (\e^{1/2})$.
The parametrization
\[
\tilde{\bx}: (x,y) \longrightarrow (x,y,u(x,y) + w_m(x,y)),
\]
valid in $B_\rho \setminus B_{c\e^{3/4}}$ becomes
\[
\tilde{\bx} : (t, \theta ) \longrightarrow (e^{-t}\cos\theta,e^{-t} 
\sin\theta, u(t,\theta) + w_m(t,\theta)) 
\]
for $(t,\theta) \in  [\bar{t}-1, -\frac{3}{4}\log \e - \log c] 
\times S^1$. We have set here $u(t,\theta) = u(e^{-t}\cos\theta,e^{-t}
\sin\theta)$ and $w_m(t,\theta) = w_m (e^{-t}\cos\theta,e^{-t}\sin
\theta)$. We now define, for all $t \in [\bar{t}-1,-\frac{3}{4}\log \e - 
\log c)$
\[
{\Sigma}_{\e,{\cal A}}(t)\equiv {\Sigma}_{\e,{\cal A}}\setminus
\tilde{\bx}\left((t,-3/4\log\e-\log c]\times S^1 \right).
\]
At this point we shall rename ${\Sigma}_{\e, {\cal A}} \equiv 
{\Sigma}_{\e,{\cal A}}(t_\e)$. Next, in ${\Sigma}_{\e,{\cal A}}$, 
we choose a unit vector field $\tilde{\nu}$ which is equal to a 
unit normal vector field  $\nu$ in ${\Sigma}_{\e, {\cal A}}(\bar{t}-1)$ 
and which equals $(0,0,1)$ in $\tilde{\bx}([\bar{t},t_\e) \times S^1)$. 
We assume that $\nu \cdot \tilde{\nu} \geq 1/2$, so that 
all surfaces near to ${\Sigma}_{\e, {\cal A}}$ are parameterized by 
${\Sigma}_{\e, {\cal A}} \ni p \longrightarrow  p + w(p)\, \tilde{\nu}(p)$,
for a suitable function $w$. The  linearized mean curvature operator, 
relative to $\tilde{\nu}$, is given by
\[
{\Lambda}_{\e,{\cal A}}\equiv \Delta_{{\Sigma}_{\e,{\cal A}}}+|{\bf
A}_{{\Sigma}_{\e,{\cal A}}}|^2 \qquad\mbox{in}\qquad{\Sigma}_{\e,
{\cal A}}(\bar{t}-1), 
\]
and by 
\begin{equation}
{\Lambda}_{\e, {\cal A}} = e^{2t}\left( \del_{tt}^2 +
\del_{\theta \theta}^2\right) + \Lambda'_u + {\Lambda}'_{\e, {\cal A}},
\label{eq:7.2}
\end{equation}
in $[{\bar t}-1 , t_\e] \times S^1$. Here $\Lambda'_u$ is the operator
in (\ref{eq:1.6}), and ${\Lambda}'_{\e, {\cal A}}$ is the correction
term coming from the geometric modifications, and in particular the
extra term $w_m$ in the parametrization for $\Sigma_{\e,{\cal A}}$. 
It is a second order operator in $t$ and $\theta$, supported in
$[\bar{t} - 1,t_\e] \times S^1$, which may be
calculated by differentiating (\ref{eq:1.8}) with respect to $w$ at 
$w=w_m$. To estimate its coefficients, we first note that the 
estimates for $w_m$ from Proposition~\ref{pr:5} translate in
the $(t,\theta)$ coefficients to 
\[
w_m = {\cal O}_\kappa(\e^{3/4}e^{-t} + \e t) \qquad
\mbox{\rm and} \qquad |\nabla^j w_m| = {\cal O}_\kappa(\e^{3/4}e^{-t}
+ \e), \quad j \geq 1.
\]
(To see this, recall that $\del_t = -r\del_r$.) Now it is not
hard to check that the coefficients of $\Lambda'_{\e,{\cal A}}$
and their derivatives are estimated by 
${\cal O}_\kappa\left(\e^{3/4}e^{t} + \e e^{2t} +\e^2 e^{4t}\right)$. 

Before discussing the mapping properties of $\Lambda_{\e,{\cal A}}$,
we define the weighted spaces on which we shall let it act.
\begin{definition} 
For $k \in {\N}$, $0<\alpha<1$ and $\delta \in {\R}$, define
${\cal C}^{k,\alpha}_\delta({\Sigma}_{\e, {\cal A}})$ by
\[
\left\{ w \in {\cal C}^{k,\alpha}_{\mathrm loc} ({\Sigma}_{\e, {\cal A}}): 
||w||_{k,\alpha,\delta}\equiv \right. 
||w||_{k,\alpha \, ({\Sigma}_{\e,{\cal A}}(\bar{t}+1))} + \sup_{
\bar{t}\leq t \leq t_\e-1} e^{-\delta t}\left. \left| w \circ 
\tilde{\bx}\right|_{k,\alpha,[t, t+1]} <\infty \right\}.
\] 
\label{de:5}
\end{definition} 

We may now state the main result of this section.
\begin{proposition}
Fix $\delta \in (1,2)$ and $\kappa >0$. Then for all $\e \in (0,\e_0]$ 
there exists an operator
\[
{\Gamma}_{\e,{\cal A}}: {\cal C}^{0,\alpha}_{\delta +2}({\Sigma}_{\e,
{\cal A}})\longrightarrow {\mathcal C}^{2,\alpha}_{\delta}({\Sigma}_{\e,
{\cal A}}),
\]
such that if $f \in {\cal C}^{0,\alpha}_{\delta+2}({\Sigma}_{\e,
{\cal A}})$, then $w=\Gamma_{\e, {\cal A}} (f)$ solves
\begin{equation} 
\left\{
\begin{array}{rllll} 
       {\Lambda}_{\e, {\cal A}} w & = &  f \qquad & \mbox{in}\qquad
{\Sigma}_{\e, {\cal A}} \\[3mm]  \pi''(w\circ \tilde{\bx}) & = &
0 \qquad & \mbox{on}\qquad \{t_\e\}\times S^1 \\[3mm] w  & = &  0 \qquad
& \mbox{on}\qquad \partial {\Sigma}_{\e, {\cal A}}\setminus
\tilde{\bx} (\{t_\e\}\times S^1).
\end{array} 
\right.
\label{eq:7.3}
\end{equation} 
Furthermore, the norm of ${\Gamma}_{\e, {\cal A}}$ is bounded 
independently of $\e$, $\kappa$ and 
all ${\cal A}$ for which $||{\cal A}|| \le \kappa \e^{3/2}$. 
\label{pr:10}
\end{proposition}

{\bf Proof :} 
The idea here is to construct a parametrix
\[
G : {\cal C}^{0,\alpha}_{\delta +2}({\Sigma}_{\e, {\cal A}}) 
\longrightarrow {\mathcal C}^{2,\alpha}_{\delta}({\Sigma}_{\e,{\cal A}})
\]
whose norm is bounded by a constant $c>0$, provided $\e$ is small enough,
by joining together two local parametrices. The first is constructed
rather explicitly inside $B_\rho$, while the second, which acts
on the exterior of this ball, is a cut-off of the solution
operator for the Jacobi operator on all of $\Sigma_0$ 
(suitably translated and rotated by ${\cal A}$), which is known
to exist by the nondegeneracy of this surface. The main point
will be to show that the norm of ${\Lambda}_{\e, {\cal A}}
\circ G -I$ can be made small, which immediately implies the result. 

In this proof, $c$ will always denote a constant which does not depend 
on $\e$, $\kappa$ or ${\cal A}$, while $c_\kappa$ may depend on
$\kappa$ but is independent of $\e$.

Appendix III contains some elementary results about the
mapping properties of $e^{2t}(\del_{tt}^2 + \del_{\theta \theta}^2)$ on
weighted spaces the cylinder $[\bar{t},t_\e] \times S^1$ which
we use now. First, from Lemma~\ref{le:2} there we obtain
a right inverse for this operator with boundary conditions $w=0$
on $\bar{t} \times S^1$ and $\pi''(w) = 0$ on $t_\e \times S^1$. 
Next, it is simple to check that
\[
\|(\Lambda'_u + \Lambda_{\e,{\cal A}}')w\|_{0,\alpha,\delta +2}
\leq (ce^{-2\bar{t}} + c_\kappa \, \e^{1/2} )\, \|w\|_{2,\alpha,\delta}. 
\]
From these two facts it is elementary to deduce the existence
of a right inverse $G^{(i)}$ for $\Lambda_{\e,{\cal A}} = e^{2t}\Delta 
+ \Lambda'_u + \Lambda'_{\e,{\cal A}}$ satisfying the appropriate boundary
conditions and with norm bounded independently of $\e$, provided
$\bar{t}$ is large enough. The superscript $(i)$ here is meant to 
connote that this is the parametrix inside the ball $B_\rho$. 

Thus the operator $\lambda (\cdot  -\bar{t})  \, G^{(i)}$ is well 
defined from ${\cal C}^{0,\alpha}_{\delta +2}({\Sigma}_{\e, {\cal A}})$ 
into ${\mathcal C}^{2,\alpha}_{\delta}({\Sigma}_{\e,{\cal A}})$ and 
has norm bounded uniformly in $\e$, for $\e$ small enough. Granted 
this, we see that the problem now reduces to solving  (\ref{eq:7.3}) 
with $f$ replaced by $g \equiv f - {\Lambda}_{\e, {\cal A}} (\lambda 
\,  G^{(i)}(f))$. The key observation is that now $g$ has support 
in  ${\Sigma}_{\e, {\cal A}} (\bar{t}+1)$, and in particular has 
a norm which is 
bounded by $c \, \|f\|_{0,\alpha, \delta + 2}$ not only in  the 
space ${\cal C}^{0,\alpha}_{\delta + 2}( {\Sigma}_{\e, {\cal A}} 
(t_\e))$ but also in ${\cal C}^{0,\alpha}({\Sigma}_{\e, {\cal A}}
(t_\e))$.

To construct the other parametrix, which is an inverse for 
$\Lambda_{\e,{\cal A}}$ outside this ball, and which we shall 
denote by $G^{(o)}$, we first make the following construction.
We modify the surface $\Sigma_{\e,{\cal A}}$ to one which has 
no boundary near zero by using the parametrization $(x,y) \to 
(x,y,u(x,y) + w_m(x,y))$ and cutting off the function $w_m(x,y)$ 
in the region $\bar{t}+1 \leq t \leq \bar{t} + 2$. More specifically, 
we let $\Sigma^c_{\e,{\cal A}}$ be the surface agreeing with 
$\Sigma_{\e,{\cal A}}$ outside $B_\rho$ and which is parametrized 
inside this ball by $(t,\theta) \to (e^{-t}\cos \theta, e^{-t}\sin
\theta,u(t,\theta) + (1 - \lambda(t - 1 -\bar{t} )) w_m(t,\theta))$. 
In ${\Sigma}^c_{\e, {\cal A}}$, we still choose a unit vector field 
$\hat{\nu}$ which is equal to the unit vector field $\tilde{\nu}$ 
in ${\Sigma}_{\e, {\cal A}} (\bar{t}+1)$ and which is equal to 
$(0,0,1)$ in the region $t\geq  \bar{t}+2$. The 
bounds for the derivatives of $w_m$ show that the surfaces 
$\Sigma_{0}$ and $\Sigma^c_{\e,{\cal A}}$ are ${\cal C}^2$ close, 
and the Jacobi operator $\Lambda^c$ for 
$\Sigma^c_{\e,{\cal A}}$ differs from that for $\Sigma_0$ by terms 
of order $c_\kappa \, \e^{3/4}$. In particular, for $\e$ small 
enough, $\Lambda^c$ is also invertible from ${\cal C}^{2, \alpha}
(\Sigma^c_{\e,{\cal A}})$ into ${\cal C}^{0, \alpha}(\Sigma^c_{\e,
{\cal A}})$ (of course, respecting the Dirichlet boundary conditions 
at the boundary of $\Sigma^c_{\e,{\cal A}}$), and we let $G^{(o)}$ 
denote its inverse whose norm is bounded uniformly in $\e$. 

We would like to have some information about the behavior of $G^{(o)}(g)$ 
near $0\in {\Sigma}^c_{\e, {\cal A}}$ when $g$ has the form specified 
above and is extended by $0$. To this aim, we apply the result of 
Lemma~\ref{le:3} in Appendix III (for example with $\delta' = \delta -2$). 
We find that there exist constants $J_0(f)$ (depending linearly
on $f$) such that 
\[
|J_0(f)|  + \| G^{(o)} (g)  - J_0 (f) \|_{2, \alpha, 
\delta - 2} \leq c_{\bar{t}} \, || f||_{0, \alpha, \delta +2}.
\]

We finally define
\[
G (f) \equiv  J_0(f) + \lambda (t_\e -\cdot ) \, ( G^{(o)} (g) - J_0(f)) + 
\lambda (\cdot -\bar{t}) G^{(i)}(f) ,
\]
where $g \equiv f - {\Lambda}_{\e, {\cal A}} (\lambda \,  G^{(i)}(f))$ 
and where we are obviously setting  $G^{(i)} = 0$ for $t \leq \bar{t}$.

We also note that $\Lambda_{\e,{\cal A}} - \Lambda^c$ is an operator 
with coefficients which are ${\cal O}_\kappa(\e^{3/4}r^{-1} + \e 
r^{-2} +\e^2 r^{-4})$ in the region $\bar{t} \leq t \leq t_\e$.
Hence, it is easy to check that $\Lambda_{\e,{\cal A}} G= I+ R$ 
where $R$ is a bounded operator 
on ${\cal C}^{0,\alpha}_{\delta + 2} (\Sigma_{\e,{\cal A}})$ with norm 
bounded by $c_\kappa \, \e^{3/4}$. As noted at the beginning, this 
suffices to complete the proof. \hfill $\Box$

\ms

Following this same proof verbatim, but replacing Lemma~\ref{le:2}
from Appendix III by Lemma~\ref{le:2.2} and using Lemma~\ref{le:3} 
with $\delta' =\delta -1$ instead of $\delta-2$,  we also obtain
\begin{proposition}
Fix $\delta \in (0,1)$ and $\kappa >0$. Then for all $\e \in 
(0,\e_0]$ there exists an operator 
\[
\hat{\Gamma}_{\e, {\cal A}}:{\cal C}^{0,\alpha}_{\delta +2}({\Sigma}_{\e,
{\cal A}}) \longrightarrow  
{\mathcal C}^{2,\alpha}_{\delta}({\Sigma}_{\e, {\cal A}}),
\]
such that for all $f \in {\cal C}^{0,\alpha}_{\delta+2 }({\Sigma}_{\e, 
{\cal A}})$, the function $w=\hat{\Gamma}_{\e, {\cal A}}(f)$ is a 
solution of the problem 
\begin{equation} 
\left\{ 
\begin{array}{rllll} 
{\Lambda}_{\e,{\cal A}} w & = &  f \qquad & \mbox{in}\qquad 
{\Sigma}_{\e, {\cal A}} \\[3mm] 
\pi'(w\circ \tilde{\bx}) & = &  0 \qquad & \mbox{on}\qquad 
\{t_\e\}\times S^1
\\[3mm]
w  & = &  0 \qquad & \mbox{on}\qquad \partial
{\Sigma}_{\e, {\cal A}}\setminus \tilde{\bx} (\{t_\e\}\times S^1).
\end{array} 
\right. 
\label{eq:7.3.3}
\end{equation} 
Furthermore, the norm of $\hat{\Gamma}_{t_\e, {\cal A}}$ is bounded 
independently of $\e$, $\kappa$ and 
all ${\cal A}$ for which $||{\cal A}|| \le \kappa \e^{3/2}$. 
\label{pr:10.10}
\end{proposition}

\ms

In Lemma~\ref{pr:7} of Appendix III we note the existence of the
bounded operator
\[
{\cal P}:\pi'' \left( {\cal C}^{2,\alpha}(S^1)\right)
\longrightarrow  {\mathcal C}^{2,\alpha}_{-2}((-\infty,0] \times S^1)
\]
such that, for any $\phi''\in \pi''\left( {\cal C}^{2,\alpha}(S^1)\right)$, 
$w={\cal P}(\phi'')$ is the unique solution in ${\cal C}^{2,\alpha}_{-2}
((-\infty, 0] \times S^1)$ of the problem
\[
\left\{ 
\begin{array}{rllll} 
\Delta w & = & 0 \qquad & \mbox{in}\quad (-\infty , 0) \times S^1\\[3mm] 
 w & =  &  \phi''   \qquad  & \mbox{on}\quad \{ 0 \} \times S^1.
\end{array}  
\right. 
\]
Now define
\begin{equation}
{\Pi}^0_{\e,{\cal A}}(\phi'')\circ \tilde{\bf x}(t,\theta) \equiv  
\lambda (t -\bar{t})\, {\cal P}(\phi'')(t - t_\e,\theta) \qquad 
\mbox{in}\qquad  [\bar{t}, t_\e] \times S^1,
\label{eq:7.33}
\end{equation}
and ${\Pi}^0_{\e,{\cal A}}(\phi'')=0$ in ${\Sigma}_{\e,{\cal A}}(\bar{t})$. 

The counterpart of Proposition~\ref{pr:8}~ is
\begin{proposition}
Fix $\delta \in (1,2)$ and $\kappa >0$. Then there exists an operator 
\[
{\Pi}_{\e,{\cal A}}:\pi''\left({\cal C}^{2,\alpha}(S^1)\right) 
\longrightarrow  
{\mathcal C}^{2,\alpha}_{\delta}({\Sigma}_{\e,{\cal A}}),
\]
such that $w= {\Pi}_{\e, {\cal A}}(\phi'')$ satisfies
\begin{equation} 
\left\{ 
\begin{array}{rllll} 
{\Lambda}_{\e, {\cal A}} w & = & 0 \qquad 
& \mbox{in}\quad 
\Sigma_{\e,{\cal A}} \\[3mm] 
\pi'' ( w \circ \tilde{\bx}) & =  & \phi'' \qquad & \mbox{on}\quad \{t_\e\}
\times S^1\\[3mm]
w  & =  &  0  \qquad & \mbox{on}\quad \partial {\Sigma}_{\e, {\cal A}}
\setminus \tilde{\bx} (\{t_\e\} \times S^1).
\end{array}  
\right. 
\label{eq:7.4}
\end{equation} 
Furthermore, $||({\Pi}_{\e,{\cal A}} - {\Pi}^0_{\e,{\cal A}}) 
(\phi'')||_{2, \alpha,\delta} \leq c_\kappa \, (\e^{3/2} + 
\e^{(3\delta +2)/4})\, ||\phi''||_{2,\alpha}$.
\label{pr:11}
\end{proposition}

{\bf Proof:} 
For simplicity, set $\tilde{w}(t,\theta) =\lambda (t -\bar{t})\,
{\cal P}(\phi'')(t -t_\e, \theta)$. The solution $w = {\Pi}_{\e,{\cal A}}
(\phi'')$ is clearly given by $\tilde{w} -\Gamma_{\e,{\cal A}}\Lambda_{\e,
{\cal A}}\tilde{w}$. It remains to estimate 
\[
||w - \tilde{w}||_{2,\alpha,\delta} \leq c_\kappa ||\Lambda_{\e,{\cal A}}
\tilde{w}||_{0,\alpha,\delta + 2}.
\]
For this we write $\tilde{w} = h(t- t_\e,\theta)$ in $[\bar{t},t_\e]
\times S^1$ and use 
\[
\Lambda_{\e,{\cal A}}\tilde{w} = e^{2 t}\Delta \lambda {\cal P}(\phi'')
(t -t_\e, \theta)+ 2 e^{2t} \nabla \lambda \cdot \nabla {\cal P}(\phi'' )
(t -t_\e, \theta) + ({\Lambda}'_u + {\Lambda}'_{\e, {\cal A}})\tilde{w}.
\]
Now, replacing $t-t_\e$ by $s \leq 0$, we see that  
\[
||\Lambda'_u \tilde{w}||_{0,\alpha,\delta+2} \leq 
c_\kappa e^{-(\delta+2)t_\e}||\tilde{\Lambda}_u' h(s)||_{0,\alpha,\delta+2}
\leq c_\kappa \e^{3(\delta + 2)/4}||\phi''||_{2,\alpha},
\]
where $\tilde{\Lambda}_u'$ is the shift by $t_\e$ of $\Lambda_u'$,
and similarly, $||\Lambda_{\e,{\cal A}}\tilde{w}||_{0,\alpha,\delta+2}$
is estimated by the same quantity. Finally, the other two terms
may be seen to be dominated by $c_\kappa \e^{3/2}||\phi''||_{2,\alpha}$
because $|h(s)| \leq c\,  e^{-2t_\e}||\phi''||_{2,\alpha} = 
c \, \e^{3/2}||\phi''||_{2,\alpha}$. \hfill $\Box$

\ms

\section{CMC surfaces near ${\Sigma}_{\e, {\cal A}}$}
\label{sect:10}

We maintain the notations of the last section. The surface parameterized by 
\[
{\Sigma}_{\e, {\cal A}} \ni p \longrightarrow p + w (p) \, \tilde{\nu} (p),
\]
has mean curvature
\begin{equation}
H = {H}_{\e, {\cal A}}+ {\Lambda}_{\e, {\cal A}} w - Q_{\e, {\cal A}}(w),
\label{eq:8.1}
\end{equation}
where ${H}_{\e,{\cal A}}$ is the mean curvature of ${\Sigma}_{\e,{\cal A}}$ 
and where $ {\cal Q}_{\e, {\cal A}}(w)$ collects the nonlinear terms. 
The form of this nonlinear term near the origin is slightly different 
than before. Indeed, the uniformity of the coefficients in (\ref{eq:1.8}) 
specifically uses the fact that the expansion for $u$ does not
have an $\e$ dependence. Thus we may not simply replace $u$ by $u + w_m$
there. Instead, we must replace $w$ by $w_m + w$ and then
expand the terms around $w_m$, with the constant and linear
terms in $w$ contributing to $H_{\e,{\cal A}}$ and $\Lambda_{\e,{\cal A}}$,
respectively. One of the terms in the expansion of $Q_u''$ about $w_m$
is a quadratic term in $w$ with a coefficient of the form 
${\cal O}(e^{2t}(|\nabla w_m| + |\nabla^2 w_m|))$. Since
$(|\nabla w_m| + |\nabla^2 w_m|) = {\cal O}(\e^{3/4} e^{-t}
+ \e)$, we see that this coefficient is ${\cal O}(\e^{3/4}e^t
+ \e e^{2t}) = {\cal O}(1+\e e^{2t})$. Hence altogether, 
\[
{\cal Q}_{\e, {\cal A}} (w) \equiv (1+ \e \, e^{2t}) \, {Q}'_{\e, {\cal
A}} ( e^t  \nabla w , e^t \nabla^2 w) +  e^t \, 
{Q}''_{\e, {\cal A}} ( e^t  \nabla w , e^t \nabla^2 w),
\]
where ${Q}'_{\e, {\cal A}}$ and ${Q}''_{\e, {\cal A}}$ are 
quadratically and cubically vanishing functions with coefficients
bounded in ${\cal C}^{k}([\bar{t},t_\e] \times S^1)$, for all $k \geq 0$, 
independently of $\kappa$, ${\cal A}$ and $\e$.

Given $\phi'' \in  \pi''\left({\cal C}^{2, \alpha}(S^1)\right) $, we wish
to construct a CMC surface which is a graph over $\Sigma_{\e,{\cal A}}$
and which has projected boundary values $\phi''$ on  
$\tilde{\bx}(\{ t_\e \} \times S^1)$. This is equivalent to solving
the boundary value problem
\begin{equation}
\left\{ \begin{array}{rlll}
{\Lambda}_{\e, {\cal A}} w & = & H_0 - {H}_{\e, {\cal A}} +  Q_{\e, 
{\cal A}} (w) \qquad & \mbox{in} \qquad {\Sigma}_{\e, {\cal A}}\\[3mm]
\pi'' ((u +w_m + w)\circ \tilde{\bf x})  & = & \phi'' \qquad & \mbox{on} \qquad 
\{ t_\e \} \times S^1 \\[3mm]
w  & = & 0  \qquad & \mbox{on} \qquad \partial {\Sigma}_{\e, {\cal A}}
\setminus \tilde{\bx} (\{t_\e\} \times S^1).
\end{array}
\right.
\label{eq:8.2}
\end{equation}

Because we are using the modified normal vector field $\tilde{\nu}$,
the surfaces we obtain will all have boundary which are vertical graphs 
over a fixed circle. Moreover, by our choice of $t_\e$ this circle 
is precisely the same one as we used for the catenoid.

Fixing $\kappa >0$, then for all $\phi'' \in \pi''({\cal C}^{2,\alpha}
(S^1))$ with $\|\phi'' \|_{2, \alpha}\leq \kappa \, \e^{3/2}$, 
we define an approximation $\tilde{w}$ to the solution of (\ref{eq:8.2}) by
\begin{equation}
\tilde{w}={\Pi}_{\e,{\cal A}}(\phi''-\pi''(\bar{w}_m(t_\e,\cdot))+
\hat{\Gamma}_{\e,{\cal A}}(H_0 - {H}_{\e,{\cal A}}),
\label{eq:8.3}
\end{equation}
which is just a solution to (\ref{eq:8.2}) if the nonlinear term is 
set to zero. We are using the function $\bar{w}_m$ from 
Proposition~\ref{pr:5.5} 
which satisfies, in particular, that $\pi'' (u+ w_m) = \pi'' \bar{w}_m$ 
on $\tilde{\bx} (\{t_\e\}\times S^1)$, but has the advantage that it
is much smaller than $u + w_m$. We are also using the right inverse 
$\hat{\Gamma}_{\e,{\cal A}}$ from Proposition~\ref{pr:10.10} here in the 
final term rather than the one from Proposition~\ref{pr:10}, which might 
be expected, simply because it affords us a better estimate, as we shall 
explain momentarily. 

Before going on, we shall collect some estimates of $\tilde{w}$.
Fix $\delta \in (1,2)$ as usual. First, let
\[
\tilde{w}_0  =  {\Pi}_{\e,{\cal A}}^0(\phi''- \pi''(\bar{w}_m(t_\e,\cdot))).
\]
We obtain from Proposition~\ref{pr:11} that 
\begin{equation}
|| \tilde{w}_0||_{2, \alpha, 2}  \leq   c \, \e^{3/2}  \, ||\phi'' 
-\pi''(\bar{w}_m ) ||_{2,\alpha}
\label{eq:8.4}
\end{equation}
and also
\begin{equation}
|| {\Pi}_{\e,{\cal A}} (\phi'' -\pi''(\bar{w}_m ) ) - \tilde{w}_0||_{2,
\alpha, \delta} \leq  c_\kappa \, (\e^{3/2} + \e^{(3\delta+2)/4}) \, 
||\phi'' -\pi''(\bar{w}_m ) ||_{2,\alpha}.
\label{eq:8.44}
\end{equation}
Furthermore, from (\ref{eq:4.5}) in  Proposition~\ref{pr:5} we get
\begin{equation}
||\pi''(\bar{w}_m)||_{2,\alpha}\leq c \, \e^{3/2}.
\label{eq:8.5}
\end{equation}
Finally, the mean curvature ${H}_{\e, {\cal A}}$ is estimated in 
Corollary~\ref{co:estmcea}. Use this estimate and also applying
Proposition~\ref{pr:10.10} with $\delta =2/3$, we have 
\begin{equation}
\| {\Gamma}'_{t_\e, {\cal A}} ( H_0 - {H}_{\e, {\cal A}} ) \|_{2, \alpha, 
2/3} \leq c \, \e^2,
\label{eq:8.6}
\end{equation}
for some constant $c >0$ which does not depend on $\kappa$. 

Putting all of these estimates together, we obtain finally that
\begin{equation}
||\tilde{w}||_{2,\alpha,[t,t+1]} \leq c\,\left( \e^2 e^{2t/3} + 
(\e^3 + \e^{(3\delta + 8)/4})e^{\delta t} + \e^3 e^{2t}\right).
\label{eq:8.66}
\end{equation}

The main reason we have had to use $\hat{\Gamma}_{\e,{\cal A}}$
rather than $\Gamma_{\e,{\cal A}}$ in (\ref{eq:8.6}) is that
otherwise the first term on the right in (\ref{eq:8.66})
would have a worse exponent, and this would lead to a far 
worse estimate in the next proposition.

Now let us solve (\ref{eq:8.2}). If we set $w = \tilde{w} + v$, then
we must prove the existence of some $v \in 
{\mathcal C}^{2,\alpha }_{\delta}({\Sigma}_{\e, {\cal A}})$ such that 
\begin{equation}
\left\{ \begin{array}{rlll}
{\Lambda}_{\e, {\cal A}} v & = &   {\cal Q}_{\e, {\cal A}} (\tilde{w} + v) 
\qquad & \mbox{in} \qquad {\Sigma}_{\e, {\cal A}} \\[3mm]
\pi'' (v \circ \tilde{\bx})  & = & 0 \qquad & \mbox{on} \qquad \{ t_\e \} 
\times S^1 \\[3mm]
v  & = & 0  \qquad & \mbox{on} \qquad \partial {\Sigma}_{\e, {\cal A}}
\setminus \tilde{\bx} (\{t_\e\} \times S^1).
\end{array}
\right.
\label{eq:8.7}
\end{equation}
As before, it is enough to find a fixed point of the mapping
\begin{equation}
{\cal M}_{\e, {\cal A}} (v) ={\Gamma}_{\e,{\cal A}}({\cal Q}_{\e, 
{\cal A}}(\tilde{w} +  v)).
\label{eq:8.8}
\end{equation}
\begin{proposition} 
For any $\kappa >0$, there exist $c_\kappa >0$ and  $\e_0 >0$ such that 
if $\e \in (0,\e_0]$ and $||\phi''||_{2,\alpha} \leq \kappa \, \e^{3/2}$,
then 
\[
||{\cal M}_{\e, {\cal A}} (0)||_{2,\alpha, \delta}\leq   c_\kappa \,  
\e^{(10+ 3\delta)/4},
\]
and
\[
||{\cal M}_{\e, {\cal A}} (v_2)- {\cal M}_{\e, {\cal A}} (v_1) 
||_{2,\alpha, \delta}
\leq \frac{1}{2} ||v_2-v_1||_{2,\alpha, \delta },
\]
provided $v_1$ and $v_2$ belong to $B \equiv \{v \, : \, ||v||_{2,\alpha , 
\delta}\leq c_\kappa \,\e^{(10+3\delta)/4} \}$. In particular, the mapping 
${\cal M}_{\e, {\cal A}}$ is a contraction on the ball $B$ into 
itself and thus ${\cal M}_{\e, {\cal A}}$ has a unique fixed point $v$ in 
this ball.
\label{pr:12}
\end{proposition}
{\bf Proof:} 
The proof is nearly identical to the proof of Proposition~\ref{pr:9}. 
The first thing we must establish is that $e^t \nabla^j \tilde{w}$ 
is bounded, and small, so that we may estimate $Q'_{\e,{\cal A}}
(e^t\nabla\tilde{w},e^t\nabla^2\tilde{w})$ by $e^{2t}(|\nabla \tilde{w}|^2
+ |\nabla^2\tilde{w}|^2)$, for example, and similarly for
the other nonlinear term. If we call the function
of $t$ on the right side of (\ref{eq:8.66}) $h(t)$, say, then
we observe that it is convex, and 
\[
h(\bar{t}) \leq c \e^2, \qquad h(t_\e) \leq c \e^{3/2}.
\]
Hence $e^t|\nabla^j\tilde{w}| \leq \e^{3/4}$, $j = 1,2$, as desired.

Now, 
\[
||{\cal M}_{\e, {\cal A}} (0)||_{2,\alpha, \delta} 
\leq c ||Q_{\e,{\cal A}}(\tilde{w})||_{0,\alpha,\delta + 2}
\]
which is estimated by the supremum of
\[
e^{-(\delta + 2)t}\left( (1+\e e^{2t})e^{2t}h(t)^2 + 
e^{4t} h(t)^3 \right) = h(t)^2 e^{-\delta t}
\left( (1 + \e e^{2t}) + e^{2t} h(t) \right).
\]
Checking the values at $t = \bar{t}$ and $t= t_\e$ and using that
the value at $\bar{t}$ also dominates the behaviour
in all of $\Sigma_{\e,{\cal A}}(\bar{t})$, we see that
\[
||{\cal M}_{\e, {\cal A}} (0)||_{2,\alpha, \delta} \leq  c \, 
\e^{(10+3\delta)/4}.
\]
This completes the proof of the first estimate.  The second
one is similar and left to the reader. \hfill $\Box$

\ms

We conclude this section with the counterpart of Corollary~\ref{co:1}. 
As we have already mentioned, we have defined $t_\e$ in such a way that
the (Dirichlet) boundary data of the surfaces defined by 
Proposition~\ref{pr:9} and Proposition~\ref{pr:13} are curves 
{\em on the same cylinder}. In the next section we shall
compare the Neumann data of the solutions of (\ref{eq:6.8}) and 
(\ref{eq:8.2}), and naturally we must differentiate
with respect to the same normal. To this aim, we note that
the relationship between the $s$ and $t$ variables on the
catenoid and surface $\Sigma_0$ is given by $ e^{-t} =\e \cosh s$
(where we assume that $t$ is close to $t_\e$ and $s$ is close to
$s_\e$). Differentiating this at $t=t_\e$, $s= s_\e$ gives 
$(dt/ds)(s_\e) = \tanh s_\e$. Since $s_\e = -(1/4)\log \e$, 
\[
\tanh s_\e \equiv \eta_\e =  \frac{1-\e^{1/2}}{1+\e^{1/2}}.
\]
We also recall the function $w_m^0$ from Proposition~\ref{pr:5.5};
in terms of the $(t,\theta)$ coordinates, 
\begin{equation}
w_m^0 (t, \theta) \equiv  e \, t + T_3 +  e^{-t} \, (R_1\, \cos \theta 
+ R_2 \, \sin \theta ) + \e \, e^t \, ( T_1\, \cos \theta +  T_2 \, 
\sin \theta ).
\label{eq:8.9}
\end{equation}
Recalling also the neighbourhood ${\cal U}$ where the parameters
${\cal A}$ reside, we set 
\[
{\cal F} \equiv {\cal U} \times \pi''\left( {\cal C}^{2,\alpha}(S^1) \right),
\]
endowed with the norm
\[
\| ({\cal A}, w ) \|_{\cal F} \equiv  || {\cal A} || + ||w||_{2, \alpha}.
\]

We now define the (slightly modified) Cauchy data mappings 
${\cal T}_\e$ for the CMC problem over $\Sigma_{\e,{\cal A}}$ and 
${\cal T}_0$ for the Laplacian on the half-cylinder $(-\infty, t_\e) 
\times S^1$:
\begin{definition}
For $\phi'' \in \pi'' \left( {\cal C}^{2,\alpha}(S^1)\right) $ with 
$||\phi''||_{2,\alpha}\leq \kappa\,\e^{3/2}$, let $w = \tilde{w} + v$ 
be the solution of (\ref{eq:8.7}) given by Proposition~\ref{pr:12}. 
Then we define 
\[
{\cal T}_\e : {\cal F} \longrightarrow  
{\cal C}^{2,\alpha} (S^1) \times  {\cal C}^{1,\alpha} (S^1) 
\]
\[
({\cal A} , \phi'') \longmapsto (\e t_\e + w_m^0 (t_\e, .)  + 
( \bar{w}_m + w)(t_\e, .), -\eta_\e \, (\e+\partial_t w_m^0 (t_\e, .) + 
\partial_t (\bar{w}_m + w )(t_\e, .))) 
\]
and 
\[
{\cal T}_0  : {\cal F} \longrightarrow 
{\cal C}^{2,\alpha} (S^1) \times  {\cal C}^{1,\alpha} (S^1)  
\]
\[
({\cal A} , \phi'') \longmapsto  ( \e t_\e  +  w^0_m (t_\e, .) + \phi'', 
- \eta_\e \, (\e  + \partial_t w^0_m + |D_\theta| \phi'' )).
\]
\label{de:7}
\end{definition}
We have made two modifications which are worth pointing out. 
First, the factor $\eta_\e$ is included so as to correspond
with differentiation with respect to $s$ on the catenoid.
Second, this is the Cauchy data with respect to the inward
pointing normal, because we are using the outward pointing
normal on the catenoid.

\begin{corollary}
For any $\kappa > 0$ there exists an $\e_0 >0$ and a constant
$c > 0$ independent of $\kappa$ such that if $\e \in (0, \e_0]$,
then ${\cal T}_\e $ and ${\cal T}_0 $ are continuous and satisfy
\begin{equation}
|| ({\cal T}_\e - {\cal T}_0)( {\cal A} ,  \phi'' ) ||_{{\cal C}^{2,\alpha} 
\times  {\cal C}^{1,\alpha}} \leq c \,  \e^{3/2}.
\label{eq:8.10}
\end{equation}
\label{co:2}
\end{corollary}
{\bf Proof}: The proof is essentially identical to the one for 
Corollary~\ref{co:1}. Continuity of the operators is
obvious. We decompose 
\[
w_m + w = \e t + w_m^0 + \bar{w}_m + \tilde{w} + v
\]
\[
= \e t + w_m^0 + \bar{w}_m + 
\Pi_{\e,{\cal A}}(\phi'') - \Pi_{\e,{\cal A}}(\pi''(\bar{w}_m
(t_\e,\cdot))) + \hat{\Gamma}_{\e,{\cal A}}(1 - H_{\e,{\cal A}}) + v.
\]
The (cut off) harmonic function on the cylinder for which ${\cal T}_0$ is
the Cauchy data operator is
\[
\e t + w_m^0 + \Pi^0_{\e,{\cal A}}(\phi'').
\]
Hence
\[
({\cal T}_\e - {\cal T}_0)({\cal A},\phi'') = 
\bar{w}_m + (\Pi_{\e,{\cal A}} - \Pi^0_{\e,{\cal A}})(\phi'')
- \Pi_{\e,{\cal A}}(\pi''(\bar{w}_m(t_\e,\cdot)))
+ \hat{\Gamma}_{\e,{\cal A}}(1 - H_{\e,{\cal A}}) + v.
\]
We estimate these in turn using Propositions~\ref{pr:5.5}, 
\ref{pr:11}, equation (\ref{eq:8.6}) and finally Proposition~\ref{pr:12}
to obtain the final estimate. \hfill $\Box$

\section{Application to our problem}

Let us now return to our original geometric problem. We are given
two CMC surfaces $\Sigma_1$ and $\Sigma_2$ which satisfy the 
assumptions of Theorem~\ref{th:1}. In this section we outline
the (very) minor changes that are needed to apply the preceding results
in our context.

First, the results of section 5 may be applied directly to the 
truncated rescaled catenoid. Similarly, we may directly apply the 
results of sections 7, 8, and 10 to the surface $\Sigma_1$. 
In particular, we obtain the corresponding mappings ${\cal T}_\e$ 
and ${\cal T}_0$, which we shall denote  
by ${\cal T}_\e^-$ and ${\cal T}_0^-$, respectively. This
superscript is meant to imply that $\Sigma_1$ is the surface
lying `underneath' $\Sigma_2$, and that its oriented normal at
the origin is $(0,0,1)$. 

However, $\Sigma_2$ is oriented oppositely, so that its
normal at the origin is $(0,0,-1)$. Thus, in section 4, 
the vector field $\bar{\nu}$ now should equal $(0,0,-1)$ in 
$\bx (B_\rho)$. The analytic modification, by adding $\e$ times
the Green function on $\Sigma_2$, and then translating vertically
by $\e a_0$, proceeds exactly as before.  The geometric
modifications of section 8 also proceed as before. However,
recall from section 3 that we had translated the catenoid
vertically by the amount $\e \log (2/\e)$, so that
its match with $\Sigma_1$ would be optimal. To make such
a match with $\Sigma_2$ at its upper boundary, we can not, of course, 
translate the catenoid again, so instead 
we translate $\Sigma_2$ vertically by the amount
$V_\e = 2\e \log (2/\e)$. The result is that 
the analogues of (\ref{eq:4.3}), (\ref{eq:4.4}) and (\ref{eq:4.5}) 
are
\begin{equation}
(x,y) \longrightarrow \left(x,y, V_\e + u_2 (x,y) - w_m (x,y) \right),
\end{equation}
\begin{equation}
(x,y) \longrightarrow \left(x,y, V_\e + \e \, \log r - \hat{w}_m (x,y) \right),
\end{equation}
and
\begin{equation}
(x,y) \longrightarrow \left (x , y, V_\e + \e \, \log r -
w^0_m - \bar{w}_m (x,y) \right),
\end{equation}
respectively, where $u_2$ is the graph function for $\Sigma_2$ and
where the functions $w_m$, $\hat{w}_m$, $w_m^0$ and $\bar{w}_m$ 
are the direct analogues of the corresponding functions for
$\Sigma_1$. We shall let the functions $w_m^0$ corresponding
to the two surfaces be denoted $(w_m^0)_\pm$, respectively. The
other functions will not need to be so explicitly labeled.

The vector field $\tilde{\nu}$ in section 9 now equals 
$(0,0,-1)$ in $\tilde{\bx} ([\bar{t}, t_\e] \times S^1)$, but
this section remains unchanged otherwise. Finally, in section 10,
the Cauchy data operators become
\[
{\cal T}^+_\e : {\cal F} \longrightarrow {\cal C}^{2,\alpha}(S^1)\times 
{\cal C}^{1,\alpha} (S^1)
\]
\[
({\cal A},\phi'') \longmapsto  (V_\e - ((w_m^0)_+ + w_m + w)(t_\e ,\cdot),
\eta_\e \, (\e + \partial_t (w_m^0)_+ (t_\e, \cdot )+ \partial_t (w_m + w)
(t_\e, \cdot)))
\]
and 
\[
{\cal T}^+_0  : {\cal F} \longrightarrow {\cal C}^{2,\alpha} (S^1) 
\times  {\cal C}^{1,\alpha} (S^1)
\]
\[
({\cal A},\phi'') \longmapsto  ( V_\e  - (w^0_m)_+ (t_\e, \cdot) -  
\phi'', \eta_\e \, (\e  + \partial_t (w^0_m)_+ (t_\e, \cdot) + 
|D_\theta| \phi'' )).
\]

\section{Matching the Cauchy data} 

We will denote by ${\cal B}'_{\kappa}$ and ${\cal B}''_{\kappa}$
the balls of radius $\kappa \, \e^{3/2}$ in the parameter space 
${\cal U}$ for ${\cal A}$ and in $\pi'' \left( {\cal C}^{2,\alpha}
(S^1)\right)$, respectively. The product ${\cal B}'_{\kappa} \times 
{\cal B}''_{\kappa}$ will be denoted simply ${\cal B}_\kappa$.
All of the constructions in the previous sections are valid for
$({\cal A},\phi'') \equiv ({\cal A}_\pm , \phi''_\pm) 
\in {\cal B}_{\kappa}^2$ for any fixed $\kappa > 0$, 
provided $\e$ is sufficiently small.

We now define the difference of the Cauchy data operators:
\[
{\bf C}_\e : {\cal B}_\kappa^2 \longrightarrow \left( {\cal C}^{2,\alpha} 
(S^1) \times {\cal C}^{1,\alpha}(S^1)\right)^2 
\]
\[
({\cal A},\phi'') \longmapsto \left( ({\cal T}_\e^+({\cal A}_+,\phi''_+) - 
{\cal S}_\e (\phi''_\pm)_+), ({\cal T}_\e^-({\cal A}_-,\phi''_-) - 
{\cal S}_\e (\phi''_\pm)_-)\right),
\]
where we have denoted by ${\cal S}_\e (\phi''_\pm)_\pm$ the 
component of ${\cal S}_\e (\phi''_\pm)$ at the upper and
lower boundaries, respectively. Setting
\[
{\cal E} = \mbox{Span}\{1, \cos \theta, \sin \theta\},
\]
then by construction, 
\[
\mbox{range }\,{\bf C}_\e \subset 
\left( {\cal E} \times {\cal C}^{1,\alpha} (S^1) \right)^2.
\]

\begin{proposition}
There exists a $\kappa_0 >0$ such that if $\kappa >\kappa_0$ then
there is an $\e_0>0$ for which, if $0 < \e < \e_0$, then
${\bf C}_\e$ has a unique zero in ${\cal B}_{\kappa}^2$.
\label{pr:13}
\end{proposition}
This Proposition produces a CMC surface $S_\e$ for each admissible 
$\e$. Indeed, if ${\bf C}_\e({\cal A},\phi'')=0$, 
then there are smooth CMC surfaces $\Sigma_1({\cal A}_-,\phi''_-)$, 
$\Sigma_2({\cal A}_+,\phi''_+)$ and the CMC perturbation of the 
truncated rescaled catenoid which we denote by 
${\cal C}_\e^c (\phi''_\pm))$, the 
union of which match up to be ${\cal C}^1$ across the two curves. 
Because of the elliptic nature of the CMC equation, it is standard that 
this union is actually ${\cal C}^\infty$ across these curves, and hence
$S_\e$ is a regular CMC surface.

Thus, in order to complete the proof of
Theorem \ref{th:1}, it remains to prove the Proposition\ref{pr:13}.

\noindent
{\bf Proof :} Let us set
\[
{\bf C}_0 :  {\cal B}_\kappa^2 \longrightarrow \left( 
{\cal C}^{2,\alpha} (S^1) \times {\cal C}^{1,\alpha} (S^1)\right)^2,
\]
\[
({\cal A},\phi'') \longmapsto \left(({\cal T}_0^+({\cal A}_+,\phi''_+)
- {\cal S}_0(\phi_\pm)_+), ({\cal T}_0^-({\cal A}_-,\phi''_-)
- {\cal S}_0(\phi_\pm)_-)\right).
\]
From Corollaries~\ref{co:1} and \ref{co:2}, we obtain
\[
\| ({\bf C}_\e -{\bf C}_0)  ( {\cal A},\phi'')\|_{({\cal C}^{2,
\alpha}\times {\cal C}^{1, \alpha})^2} \leq c\, \e^{3/2},
\] 
where the constant $c >0$ does not depend on $\kappa$. 

We examine the map ${\bf C}_0$ more closely. Recall first the small
deviation, of order ${\cal O}(\e r)$ of the translated 
catenoid from $-\e \log r$. This error term is clearly radial,
since the catenoid is rotationally symmetric, and hence
we write it as $\e \beta(r)$, where $\beta(0) = 0$. Now
\[
\begin{array}{rlll}
{\bf C}_0({\cal A},\phi'') & = &  \left( ((w_m^0)_+ + \e \beta(r),
\eta_\e(\del_t( (w_m^0)_+ + \e \beta(r))) - (\eta_\e - 1)
|D_\theta|\phi''_+), \right. \\[3mm]
&  & \left. ((w_m^0)_- - \e \beta(r), \eta_\e(\del_t( (w_m^0)_- 
+ \e \beta(r))) + (\eta_\e - 1)|D_\theta| \phi''_-)\right).
\end{array}
\]
It is trivial to see that ${\bf C}_0$ is an isomorphism from 
$({\mathbb R}^6 \times \pi'' ({\cal C}^{2, \alpha}(S^1)))^2$ into 
$({\cal E}\times {\cal C}^{1, \alpha}(S^1))^2$.
In particular, there is a unique zero of this mapping,
namely where $\phi'' = 0$ and the $(w_m^0)_\pm$ are chosen to
cancel $\pm \e \beta(r)$. Notice that this solution is certainly
within ${\cal B}_\kappa^2$, because $|\e \beta(r)| \sim \e^{7/4}$.

We would like to use a degree theoretic argument to conclude
that there is also a single zero of ${\bf C}_\e$ within
${\cal B}_\kappa^2$. Unfortunately, the nonlinear correction
terms in the difference ${\bf C}_\e - {\bf C}_0$, whilst
small, are not compact. We could, of course, use a contraction
mapping argument again, but we propose, instead, the following
shorter route. We write
\[
\begin{array}{rllll}
{\bf C}_\e({\cal A},\phi'') & = & \left(w_m^0 + \e \beta(r) + 
F'({\cal A},\phi'') + F''({\cal A},\phi''), \right. \\[3mm]
& &  \left. 
\del_s(w_m^0 + \e \beta(r)) + (\eta_\e - 1)|D_\theta|\phi''
+ \del_s F'({\cal A},\phi'') + \del_s F''({\cal A},\phi'')\right).
\end{array}
\]
Here $F' = (I - \pi'')(({\bf C}_\e - {\bf C}_0)({\cal A},\phi''))$,
and $F'' = \pi''(({\bf C}_\e - {\bf C}_0)({\cal A},\phi''))$. 
The range of $F'$ lies in the
finite dimensional space ${\cal E}^2$, but the range of $F''$
is ostensibly the problem, since it is infinite dimensional.
These error terms are all, however, ${\cal O}(\e^{3/2})$,
with constants independent of $\kappa$.

Define a family of smoothings of this map ${\bf C}_{\e,q}$,
for $0 < q < 1$, by replacing the terms
$F''({\cal A},\phi'')$ and $\del_s F''({\cal A},\phi'')$
by $|D_\theta|^{-q}F''({\cal A},\phi'')$ and 
$|D_\theta|^{-q}\del_s F''({\cal A},\phi'')$, respectively. 
Here $|D_\theta|^{-q}$ is the pseudodifferential operator
of order $-q$ defined by
\[
|D_\theta|^{-q}: \sum_{|n| \geq 2}  a_n e^{in\theta}
\longrightarrow \sum_{|n| \geq 2}  |n|^{-q} a_n e^{in\theta}.
\]
Since the norm of $|D_\theta|^{-q}$, when defined from 
${\cal C}^{1, \alpha}(S^1)$ into itself, is bounded independently 
of $q$ for $0 < q < 1$, we see that the nonlinear terms are
all still ${\cal O}(\e^{3/2})$, independently of $\kappa$.

It is now commonplace, using the Leray-Schauder degree,
that there exists  $({\cal A},\phi'')_q$ for 
which ${\bf C}_{\e,q}(({\cal A},\phi'')_q) = 0$. More 
specifically, this point exists in ${\cal B}_\kappa^2$ if we
first choose $\kappa$ large enough to overwhelm the other
($\kappa$-independent) constants which estimate the nonlinear
terms in the mapping, and then choose $\e$ accordingly 
sufficiently small. We also note that 
$\|({\cal A},\phi'')_q\|$ is always in ${\cal B}_\kappa^2$,
hence has norm bounded uniformly in $q$. This means
that we may extract a sequence $q_j \rightarrow 0$ such
that $({\cal A},\phi'')_{q_j}$ converges in 
${\cal U} \times {\cal C}^{2,\alpha'}$ for any fixed
$\alpha' < \alpha$. This is clearly sufficient for
our purposes, and it is clear that the limit of this
sequence is a zero of ${\bf C}_\e$. This completes our proof.
\hfill $\Box$

\ms

\section{Technical information needed for the proof of generic 
nondegeneracy}

Now that we have proven the existence of the family of 
CMC connected sums $S_\e$ of the two surfaces $\Sigma_1$
and $\Sigma_2$, we turn our attention to establishing criteria 
ensuring that the $S_\e$ are nondegenerate. 
This will require some preparatory work. In this section we 
give estimates on the graph function for $S_\e$ over
the truncated rescaled catenoid and use this to describe the form 
of the Jacobi operator on $S_\e$. In the next section we
give precise estimates for the solutions of this Jacobi
operator corresponding to the low eigenmodes $j = 0, \pm 1$
on the cross-section. After that we will be able to address
the nondegeneracy question directly.

In the previous sections we gave good estimates for $S_\e$
as a graph over the truncated rescaled catenoid $\Sigma_\e^c$,
specifically in the region where the parameter $s$ lies
in $[-s_\e,s_\e]$ ($s_\e = -(1/4)\log \e$).  However,
we shall need to extend these estimates to the larger
region including the balls $B_\rho \setminus 
B_{c\e^{3/4}}$ in each of the surfaces $\Sigma_j$.
This entire region may be written as a graph over 
a region in $\Sigma_\e^c$. Recall the relationships
between the various variables we have used:
\[
r = e^{-t} \qquad \mbox{\rm and } \qquad e^{-t} =
\e \cosh s.
\]
Since the annuli in $\Sigma_j$ are parametrized by
$[\bar{t},t_\e] \times S^1$, then if we define $\bar{s}_\e$ by
\[
e^{-\bar{t}} = \e \cosh \bar{s}_\e,
\]
we see that the region in $S_\e$ of interest to us, which
we write as $S_\e \cap B_\rho$, is
parametrized by $[-\bar{s}_\e,\bar{s}_\e] \times S^1$.

Notice also that $S_\e \cap B_\rho$ may be decomposed into
three components. The first central component, denoted by $I$, 
corresponds to $s$ lying in the interval $[-s_\e, s_\e]$.
The two other components, $II_1$ and $II_2$, are vertical graphs 
over $\Sigma_1$ and $\Sigma_2$, respectively.

\begin{lemma}
For some small value of $\rho$, and for $\e$ sufficiently small,
there is a function $g_\e$ on $\Sigma_\e^c$ such that
\[
{\bf x}_\e : [- \bar{s}_\e, \bar{s}_\e] \times S^1 \ni (s,\theta)
\longrightarrow {\bf x}^c_\e (s , \theta) + g_\e (s, \theta) \, 
\bar{n}_\e (s, \theta) \in S_\e,
\]
where $\bar{n}_\e$ is the unit vector field on $\Sigma_\e^c$
defined in (\ref{eq:6.3}). Furthermore, the estimate
\begin{equation}
\nabla^k g_\e (s, \theta) = {\cal O} (\e^{2} \cosh^2 s ),
\label{eq:11.1}
\end{equation}
holds for $(s, \theta) \in [ - \bar{s}_\e , \bar{s}_\e ] \times S^1$ 
when $k\geq 1$ but only for $(s, \theta) \in [-s_\e , s_\e ] \times
S^1$ when $k=0$.
\label{le:7}
\end{lemma}
{\bf Proof:} 
In the region $I$, where $|s| \leq s_\e$, $g_\e = 
{\tilde w} +v$ as in Proposition~\ref{pr:9}, and so 
(\ref{eq:11.1}) follows directly from (\ref{eq:6.10}) and 
(\ref{eq:6.11}).

In the regions $II_i$, when $s_\e \leq s \leq \bar{s}_\e$, 
$g_\e = \hat{w}_m + \tilde{w} + v$, and so 
we use the estimates in Proposition~\ref{pr:12}, (2) of 
Proposition~\ref{pr:5} and (\ref{eq:8.66}). \hfill $\Box$

The restriction to $k\geq 1$ in the outer shell is simply
because of the presence of the term $\e^{3/2}|\log \e|$ 
when $k=0$.

\ms

We also need the 
\begin{lemma}
For $i = 1,2$, the component $II_i$ can be parametrized by
\[
{\bf x}_{i,\e} : (x,y) \in B_\rho\setminus B_{\e^{3/4}} \longrightarrow
(x,y,u_i(x,y)+ h_{i,\e}(x,y)),
\]
where $h_{i,\e}$ satisfies
\begin{equation}
\nabla^k h_{i,\e} (x,y) = {\cal O} \left( r^{-k} (\e + \e^{3/4} r) \right)
\label{eq:11.2.2}
\end{equation}
for $k \geq 1$. 
\label{le:7.7}
\end{lemma}
The proof is similar to the proof of Lemma~\ref{le:7}; the only 
difference is that $\hat{w}_m$ must be replaced by $w_m$. Details
will be omitted.  Again the restriction to $k\geq 1$ is simply
to avoid a logarithmic term when $k= 0$.

Finally, recall that $- (\e^2 \cosh^2 s)^{-1} {\cal L}$ is the 
Jacobi operator about $\Sigma_\e^c$ with respect to the normal 
vector field $n$ while as in the expression following (\ref{eq:6.7}),
$(\e^2 \cosh^2 s)^{-1}\left(- {\cal L} + L_\e\right)$ 
is the Jacobi operator about this same surface with respect to the 
transverse vector field $\bar{n}_\e$. The coefficients of $L_\e$ are of 
order $1/(\cosh s)^2$ and are supported in the region $s_\e - 2 \leq |s| 
\leq s_\e - 1$. We now let $\LL_\e$ be the Jacobi operator on 
$S_\e$ with respect to $\bar{n}_\e$. 
\begin{corollary}
When $(s, \theta) \in [-\bar{s}_\e, \bar{s}_\e]\times S^1$, 
\[
{\LL}_\e = -\frac{1}{\e^2 \cosh^2 s} \left( {\cal L} - L_\e + 
{\LL}'_{\e}\right),
\]
where ${\LL}'_\e$ is a second order operator the coefficients of which, 
along with their derivatives, can be estimated by a constant times 
$(\e + \e^2 \, \cosh^2 s )$ for $(s, \theta) \in [ -\bar{s}_\e, 
\bar{s}_\e ] \times S^1$. 
\label{co:3}
\end{corollary}
{\bf Proof :} Following (\ref{eq:6.7}), the mean curvature 
of any graph over the catenoid $\Sigma^c_\e$, parametrized 
using $\bar{n}_\e$, is given by
\[
\begin{array}{rllll}
- \ds \frac{1}{\e^2 \cosh^2 s} \, & \ds \left({\cal L} w  -  L_\e w + 
\e \bar{Q}'_\e \left(\frac{w}{\e \cosh s}, \frac{\nabla w}{\e \cosh s}, 
\frac{\nabla^2 w}{\e \cosh s} \right) \right.\\[3mm]
& \ds \quad \left. + \e \cosh s \, \bar{Q}''_\e \left(\frac{w}{\e \cosh s},\frac{\nabla w}{\e \cosh s},\frac{\nabla^2 
w}{\e \cosh s}\right) \right),
\end{array}
\] 
The operator $\LL_\e'$ is obtained by linearizing the last
two expressions around $g_\e$. Notice that when $|s| \geq s_\e$, 
$\bar{n}_\e$ is identically equal to $(0,0,\pm 1)$ and 
so by (\ref{eq:1.5}), the nonlinear terms only involve
the derivatives of $g_\e$ and not $g_\e$ itself in this
range, which means that we may use the estimate (\ref{eq:11.1})
in this region.   \hfill $\Box$

\ms

\section{Jacobi fields}

As we discussed at the beginning of the last section, we
require precise asymptotics for the Jacobi fields for 
$\LL_\e$ corresponding to the low eigenmodes on the
circle. More specifically, there are explicit Jacobi fields on the 
catenoid, i.e. solutions of ${\cal L} w=0$ in ${\R}\times S^1$, 
given by 
\[ 
\Psi^{0,+} (s, \theta) = \tanh s, \qquad \qquad \qquad  \Psi^{0,-} (s, 
\theta) = (1 - s \tanh s),
\]
\[
\Psi^{\pm 1,+} (s, \theta) = \displaystyle{\frac{1}{\cosh s}} 
e^{\pm i \theta}, \qquad \mbox{and}\qquad \Psi^{\pm 1,-}(s,\theta) 
= (\displaystyle{\frac{s}{\cosh s}} + \sinh s)e^{\pm i\theta}.
\]
These all arise from explicit families of perturbations of
the catenoid. In fact, if $S$ is any CMC surface and $S(\eta)$ is a 
smooth one-parameter family of CMC deformations with $S(0) = S$, 
then $S(\eta)$ may be written as a graph (with respect to
some transverse normal vector field) over $S$ for small $\eta$.
Actually, all that is needed is that this graph function exist
over any fixed compact set of $S$ for some nontrivial range of 
values of $\eta$ which might diminish to zero as the compact set grows.
This is sufficient to make sense of the derivative of the graph
function at $\eta = 0$, and this derivative is a Jacobi field.
The Jacobi fields above are obtained in this way, as derivatives
of one parameter families of CMC surfaces parametrized using
the unit normal vector field; $\Psi^{0,+}$ and $\Psi^{\pm 1,+}$
correspond to vertical and horizontal translations, respectively, 
$\Psi^{0,-}$ corresponds to changes by dilation and $\Psi^{\pm 1,-}$
correspond to rotations about the $x$ and $y$ axes.  
If we write the graphs using the vector field $\bar{n}_\e$
instead, then the corresponding Jacobi fields will be denoted 
$\bar{\Psi}^{j,\pm}_\e$. These are solutions of $({\cal L} - L_\e)w = 0$, 
and from (\ref{eq:1.11}) in Appendix I we have 
\[
\bar{\Psi}^{j,\pm}_\e = \frac{1}{n \cdot \bar{n}_\e} \, \Psi^{j, \pm}.
\]

The goal of this section is to find good estimates for the
Jacobi fields on the surfaces $S_\e$ which are perturbations
of these; these will be solutions of $\LL_\e w = 0$ and will be
denoted by $\Phi^{j, \pm}_\e$ for $j = 0, \pm 1$.  We are really
only interested in describing them over the regions $II_i$, $i=1,2$.

The five Jacobi fields which correspond to vertical and horizontal 
translations and rotations of the vertical axis 
are the easiest to describe. We shall only need to 
describe their behaviour over the regions $II_i$, and will now use
the variables $(x,y)$ rather than $(s,\theta)$ there.

\begin{proposition} 
The Jacobi fields $\Phi^{j, +}_\e$, $j = 0, \pm 1$, and $\Phi^{\pm 1, -}_\e$ 
are described in $II_i$, $i=1,2$, by
\[
\Phi^{0,+}_\e (x,y) = (-1)^{i},
\]
\[
\Phi^{+1, +}_\e (x,y) = (-1)^i  \, \del_x \tilde{u}_{i,\e}(x,y)) ,
\qquad
\Phi^{-1, +}_\e (x,y) = (-1)^i \, \del_y \tilde{u}_{i,\e}(x,y))
\]
and 
\[
\Phi^{+1,-}_\e (x,y) = (-1)^i \,(x + \tilde{u}_{i,\e}(x,y)\,\del_x 
\tilde{u}_{i,\e}(x,y)),
\quad 
\Phi^{-1, -}_\e (x,y) = (-1)^i \, (y + \tilde{u}_{i,\e}(x,y) \, \del_y \tilde{u}_{i,\e}(x,y)).
\]
where $\tilde{u}_{i, \e}\equiv u_i +h_{i, \e}$.
\label{pr:4.11} 
\end{proposition} 
{\bf Proof:} 
The simple expression for $\Phi_\e^{0,+}$ follows from the
fact that $\bar{n}_\e = (0,0,(-1)^{i+1})$ in $II_i$. On the
other hand, recall from Lemma~\ref{le:7.7} that in these regions the 
graph functions for $S_\e$ relative to $\bar{n}_\e$ have the 
form $u_i(x,y) + h_{i,\e}(x,y)$. Differentiating
with respect to $x$ and $y$ corresponds to infinitesimal
translations in these directions, and this leads to the
stated expressions. The Jacobi fields corresponding to the two
rotations of the vertical axis can be obtained similarly.
\hfill $\Box$

\ms

Unfortunately, it is more difficult to get good estimates for
the last remaining Jacobi field since we have not proved that 
$S_\e$ depends smoothly on $\e$. We will obtain this last function,
and estimates for it, by a perturbation argument.

\begin{proposition}
Assume that  $\delta \in (1,2)$. Then for some $\bar{s}_1 >0$ sufficiently
large, but independent of $\e$, and $\e$ is small enough, there exists 
a Jacobi field $\Phi^{0,-}_{\e}$, 
defined in $[-\bar{s}_\e+ \bar{s}_1 , \bar{s}_\e - \bar{s}_1]\times S^1$, 
which satisfy
\[
\Phi^{0,-}_\e (x,y)= - \log (2r/\e) + {\cal O} (r + r^\delta \, | \log \e | ),
\]
in ${\bx}_{i,\e}(B_{\bar{\rho}_1}\setminus B_{\e^{3/4}})$, for $i=1,2$. 
By definition here, $\bar{\rho}_1 \equiv \e\, \cosh (\bar{s}_\e -\bar{s}_1)$.
\label{pr:14}
\end{proposition}
{\bf Proof:} First, by (\ref{eq:1.11}) in Appendix I, 
\begin{equation}
\displaystyle \LL_\e (\frac{1}{n\cdot \bar{n}_\e} \, w )= -\frac{1}{\e^2 \cosh^2 s} {\cal L} w + 
\frac{1}{\e^2 
\cosh^2s}{\LL}''_{\e} w
\label{eq:cvf}
\end{equation}
where the operator $\LL''_\e$ enjoys the same properties as $\LL'_\e$, namely 
has all its coefficients bounded by a constant times $\e + \e^2 \, \cosh^2 s$. 
Therefore, it is enough to find the appropriate Jacobi fields for the operator 
\[
{\cal L} - {\LL}''_{\e}.
\]

First, if $\bar{s}_1$ is chosen large enough, the result of 
Proposition~\ref{pr:6} holds for all $\e$ small enough, 
with $s_0 = \bar{s}_\e-\bar{s}_1$ and with ${\cal L}$ replaced by ${\cal L} - 
{\LL}''_{\e}$. Indeed, 
\[
\|{\LL}''_\e (w)\|_{0, \alpha , \delta} \leq c \, (\e+ e^{-2\bar{s}_1})\, 
\|w\|_{2, \alpha, \delta}.
\]
The claim follows immediately, provided $\bar{s}_1$ is chosen
large enough. From now on we keep $\bar{s}_1$ fixed so that this 
is true and we will denote by ${\mathbb G}_{\bar{s}_\e -\bar{s_1}}$ 
the right inverse obtained by perturbing the right
inverse ${\cal G}_{\bar{s}_\e  -\bar{s}_1}$ for ${\cal L}$.

We obtain the desired function
\[
\Phi^{0,-}_\e = \Psi^{0, -} - {\mathbb G}_{\bar{s}_\e-\bar{s}_1} 
({\LL}''_\e (\Psi^{0,-}))
\]
easily enough. 

The main work will be in estimating ${\cal G}_{\bar{s}_\e - \bar{s}_1} 
({\LL}''_\e (\Psi^{0,-}))$. First, recall that $r =\e \cosh s$. 
This implies that when $s>0$, 
\[
\e \, e^{s}= 2 \, r+{\cal O}(\e^2 \, r^{-1}),\qquad \e \, e^{-s}=
{\cal O}(\e^2 \, r^{-1}), \qquad \mbox{\rm and } \qquad
s = \log (2r/\e) +{\cal O}(\e^2/r^2),
\]
and so
\[
1- s \, \tanh s = 1 - \log (2 r / \e) +{\cal O}(\e^2 \, 
| \log \e | \, r^{-2}).
\]
On the other hand, when $s <0$, 
\[
\e \, e^{-s}= 2 \, r +{\cal O} ( \e^2 \, r^{-1}),\qquad \e \, e^{s}={\cal O}(\e^2 \, r^{-1}),
\qquad \mbox{\rm and } \qquad s = - \log (2r/\e) +{\cal O}(\e^2/r^2),
\]
which gives
\[
1- s \, \tanh s = 1 - \log (2 r / \e) +{\cal O}(\e^2 \, |\log \e | \, r^{-2}).
\]

We next show that we can get somewhat sharper estimates for
${\cal G}_{\bar{s}_\e -\bar{s}_1}({\LL}''_\e (\Psi^{j,\pm}))$ 
than those obtained from Proposition~\ref{pr:6} directly.  
Using the bounds on the coefficients of $\LL_\e''$ we find that
\[
|({\LL}''_\e (\Psi^{0,-})|\leq c \, (\e + \e^2 \, \cosh^2 s)\, 
(1+ |s|)\leq c \, (\e + \e^\delta \, \cosh^\delta s)\, | \log \e | ,
\]
for some constants $c>0$ which are independent of $\e$. We have
also estimated $(\e \cosh s)^k$, $k = 2, 3$, by $\e^\delta \cosh^\delta s$ 
here in order to simplify later estimates. 

Now recall the construction of Proposition~\ref{pr:6}. Let us write 
\[
w = {\cal G}_{\bar{s}_\e -\bar{s}_1} ({\LL}''_\e (\Psi^{0,-})) 
=\sum_{n\in {\Z}} w_n(s)\, e^{in\theta}, \qquad
f = {\LL}''_\e (\Psi^{0,-}) =\sum_{n\in {\Z}} f_n(s) \, e^{in\theta}.
\]
As in that proof, when $|n|\geq 2$, multiples of the function $\, n^{-2}\, 
(\e+ \e^\delta \, (\cosh s)^{\delta})\, |\log \e |$ can be used as
supersolutions for $\pm w_n$. Hence, for $|n|\geq 2$, 
\[
|w_n(s)| \leq \frac{c}{n^2} \, \left( \e + \e^\delta  \, 
\cosh^{\delta} s\right) \, |\log\e|.
\]
To handle the remaining cases $n = 0, \pm 1$ we use the explicit
formul\ae\ (\ref{eq:5.6}) and $(\ref{eq:5.7})$ 
\[
w_0(s) = \tanh s \int_0^s \tanh^{-2} t \int_0^t \tanh u \, 
f_0 (u) \, du \, dt,
\]
and 
\begin{equation}
w_{\pm 1}(s) = \cosh^{-1} s \int_0^s \cosh^2 t \int_0^t \cosh^{-1} u \, 
f_{\pm 1}(u) \, du \, dt.
\label{eq:wi}
\end{equation}

Direct estimates yield 
\[
|w_j(s)| \leq c \, \left( \e |\log \e|^3 +\e \cosh s + \e^\delta \cosh^\delta s 
|\log \e| \right) \leq c \left( \e \cosh s  + r^\delta |\log \e|\right),
\qquad j = 0, \pm 1.
\]

Summation over $n$ now yields the desired estimate for
the remainder term. The derivatives are handled similarly. \hfill $\Box$

\ms

\section{Proof of generic nondegeneracy}

Fix $(p_1,p_2 ,\theta) \in \Si_1 \times \Si_2 \times S^1$,
and then choose rigid motions of the surfaces $\Si_j$ so that the
points $p_j$ are mapped to the origin and the tangent planes
$T_{p_j}\Si_j$ are mapped to the $x\, y$-plane with opposite orientation.
Suppose furthermore that we first normalize these mappings so that
the principal directions at these points are mapped to the
$x$ and $y$ axes, respectively. (There is of course a choice
to be made here regarding the ordering of the principal directions,
but we require, for example, that the direction with larger
principal curvature be carried to the $x$-axis; since we are
specifying an orientation, this fixes the choice at all
points except umbilics.) Finally, rotate $\Si_2$ about the
$z$-axis by an angle $\theta$ so that its principal directions are 
aligned with the vectors $(\cos \theta,\sin\theta,0)$ and $(-\sin\theta,
\cos\theta,0)$. We call the resulting singular configuration
$\Si_1 \sqcup \Si_2(p_1,p_2,\theta)$. The resulting `moduli space',
$C(\Si_1, \Si_2)$, of such configurations is clearly five dimensional. 
It is the quotient of an eleven-dimensional space 
by the (six-dimensional) group of rigid motions. This procedure
yields local charts on $C(\Si_1, \Si_2)$. 
Finally, given some sufficiently small $\e > 0$, we form
the desingularized connected sum $S_\e(p_1,p_2,\theta)$. 
Note that $\Si_1 \sqcup \Si_2(p_1,p_2,\theta)$ is the union, near
the origin, of two graphs over the $x\, y$-plane
\[
{\bx}_i : B_\rho \ni (x,y) \longrightarrow (x,y, u_i (x,y)), \qquad
i = 1, 2.
\]
The maps ${\bx}_i$ and $u_i$ depend, of course, on $p_1, p_2$ and $\theta$.

Our aim in this final section is to prove Proposition~\ref{pr:2},
that is, to prove the nondegeneracy of $S_\e(p_1,p_2,\theta)$ for 
$\e$ small. Recall that this means that we need to show that
there are no nontrivial Jacobi fields on $S_\e(p_1,p_2,\theta)$
which vanish on $\partial S_\e(p_1,p_2,\theta)$.  We are not
able to show that this is true for every value of the parameters,
but at least we shall show that it holds generically, in a
precise sense.

We first prove a result which gives a criterion for 
nondegeneracy. 
\begin{theorem}
Let $(p_1,p_2,\theta) \in \Sigma_1 \times \Sigma_2 \times S^1$ be fixed. 
If there exists a sequence $\e_n \to 0$ for which the 
surface $S_{\e_n} (p_1, p_2, \theta)$ is degenerate, then
\begin{equation}
\left(\del_{xx}^2 u_2 (0,0)-\del_{xx}^2 u_1 (0,0)\right)
\left(\del_{yy}^2 u_2 (0,0)-\del_{yy}^2 u_1(0,0)\right)
- \left(\del_{xy}^2 u_2 (0,0) -\del_{xy}^2 u_1 (0,0)\right)^2 =0 .
\label{eq:nondegeqn}
\end{equation}
\label{th:last}
\end{theorem}
{\bf Proof:} We omit $p_1,p_2$ and $\theta$ from the notation
since they are fixed. Let ${\LL}_\e $ denote the mean curvature operator
linearized about $S_\e$ with respect to the normal transversal 
vector field used in the previous section. We shall also often
simply write $\e$ instead of $\e_n$. The degeneracy of $S_\e$ 
means that there exists a nontrivial function $w_\e$ on $S_\e$ 
with $w_\e =0$ on $\del S_\e$ and such that ${\LL}_\e \, w_\e = 0$. 

Fix any $\delta_0 \in (1,2)$. We now choose, for each $\e$, 
a weight function $\gamma_\e : S_\e \rightarrow {\R}$ which satisfies
\[
\gamma_\e (p) \equiv  1 \qquad \mbox{in}\qquad  S_\e \setminus B_{2\rho} (0),
\qquad \qquad \gamma_\e \circ {\bf x}_{i,\e}(t,\theta) \equiv e^{\delta_0 t}
\qquad \mbox{in} \qquad [\bar{t}_1 , t_\e] \times S^1,
\]
for $i = 1, 2$, and 
\[
\gamma_\e \circ {\bf x}_\e (s, \theta ) \equiv  (\e \, 
\cosh s)^{-\delta_0} \qquad \mbox{in} \qquad [-s_\e  , s_\e] \times S^1.
\]
We also require that $\gamma_\e$ and its derivative are bounded 
independently of $\e$ in $S_\e \cap (B_{2\rho} \setminus B_\rho )$. 

Use these weight functions to normalize the functions $w_\e$ by
\[
\sup_{p \in S_\e} \gamma_\e (p) \, |w_\e (p)| = 1.
\]
Suppose that $p_\e \in S_\e$ is a point where this supremum is 
achieved. Passing to a subsequence, we may assume that $\{p_\e\}$ 
converges to some point $p_\infty \in \Sigma_1 \cup \Sigma_2$. 
We distinguish various cases according to the location of $p_\infty$.

\noindent
{\bf Case 1.}
Assume that $p_\infty =0$. In this case, we may write (at least for 
$\e$ small enough)
\[
p_\e = {\bf x}_\e (s'_\e,\theta'_\e),
\]
for some $(s'_\e,\theta'_\e) \in  [\log \e + c , -\log \e - c ] 
\times S^1$. We distinguish two further cases according to 
the behaviour of the sequence $s'_\e$.

\noindent
{\bf Subcase 1.1.} Assume that (up to a subsequence) 
$(s'_\e,\theta'_\e)$ converges to $(s_0, \theta_0) \in 
{\R}\times S^1$. Then define
\[
\tilde{w}_\e (s,\theta) = \e^{\delta_0} \, w_\e \circ {\bf x}_\e (s, \theta).
\]
This still solves ${\LL}_\e \tilde{w}_\e =0$ in $[\log \e - c , 
- \log \e + c ] \times S^1$, is bounded by $(\cosh s)^{-\delta_0}$ 
and also satisfies
\[
\tilde{w}_\e (s'_\e,\theta'_\e)\equiv 1.
\]
Now pass to the limit, possibly after passing to a further subsequence.
By Corollary~\ref{co:3} we obtain a nontrivial function $w$ such that
\begin{equation}
\del_{ss}^2 w + \del_{\theta\theta}^2 w + \frac{2}{\cosh^2 s} w=0
\label{eq:11.5}
\end{equation}
in ${\R}\times S^1$ and which is bounded by  $(\cosh s)^{-\delta_0}$. 
We now show that this is not possible. Let 
\[
w (s , \theta) =\sum_{n \in {\Z}} w_n (s) e^{in \theta}
\]
be the Fourier decomposition of $w$. Then 
\[
\hat{w}(s , \theta) =\sum_{|n| \geq 2 } w_n (s) e^{in \theta},
\]
is still a solution of (\ref{eq:11.5}); it also decays exponentially 
at both $\pm \infty$. Multiplying (\ref{eq:11.5}) by $\hat{w}$ and 
integrating by parts we find
\[
\int_{{\R}\times S^1} \left( (\del_s \hat{w})^2 + (\del_\theta 
\hat{w})^2 - \frac{2}{\cosh^2 s} \hat{w}^2 \right) ds \, d\theta =0,
\]
which implies that $\hat{w}=0$. Hence $w = \sum_{|n| \leq 1} 
w_n(s) e^{in\theta}$. As in the last section, the solutions in
these low eigenspaces are linear combinations of the 
explicit solutions $\Psi^{j,\pm}$, $j= 0, \pm 1$, and no
nontrivial solution of this form can decay as quickly as
$(\cosh s)^{-\delta_0}$ at $\pm \infty$. Hence this subcase
cannot occur.

\noindent
{\bf Subcase 1.2.} Now assume that $\lim_{\e \rightarrow 0} 
s'_\e = +\infty$ or $-\infty$. To fix ideas, assume that $\lim_{\e 
\rightarrow 0} s'_\e = - \infty$. Notice that because $\lim_{\e \to 0}
p_\e = 0$, we also have $\lim_{\e \to 0} s'_\e - \log \e = +\infty$. 
Define
\[
\hat{w}_\e (s, \theta) = \e^{\delta_0} \, (\cosh s'_\e )^{\delta_0}\,
w_\e (s +s'_\e, \theta).
\] 
This function is bounded by a constant times $e^{\delta_0 s}$ in 
$[\log \e - s'_\e, -s'_\e]\times S^1$ and satisfies
\[
\lim_{\e \rightarrow 0} \hat{w}_\e (0, \hat{\theta}_\e) =1.
\] 
Again passing to the limit as $\e \to 0$ using Corollary~\ref{co:3}, 
we obtain a nontrivial solution of
\[
\del_{ss}^2 w + \del_{\theta\theta}^2 w =0 \qquad \mbox{in} 
\qquad {\R} \times S^1,
\]
which is bounded by $e^{-\delta_0 s}$. Since $\delta_0 \notin \Z$, 
this is impossible, which rules out this subcase.

\noindent
{\bf Case 2.} Finally we assume that $\lim_{\e \rightarrow 0} p_\e 
\neq 0$. Possibly extracting subsequences, we pass to the 
limit as $\e$ tends to $0$ and obtain two solutions $w_1$ and $w_2$ 
(at least one of which is nontrivial) of
\[
\Lambda_i w_i = 0 \qquad \mbox{in} \qquad \Sigma_i \setminus \{0\},
\]
with $w_i=0$ on $\partial \Sigma_i$. We know that $w_i \in 
{\cal C}^{2,\alpha}_{- \delta_0}(\Sigma_i \setminus \{0\})$, and
so there must exist constants $c^i_j \in {\R}$, $j = 0, \pm 1$, such that
\[
{\Lambda}_i w_i = - 2 \pi \, \left( c^i_0 \delta_0 + (c^i_{1},c^i_{2})  
\cdot \nabla \delta_0 \right).
\]
Our goal is to show that these constants $c^i_j$ all vanish. We claim
that this follows from the condition (\ref{eq:nondegeqn}). 
Granting this, then each $w_i$ must be a regular Jacobi fields over 
the whole of $\Sigma_i$, and at least one of them must be
nontrivial. Nondegeneracy of the two surfaces implies 
that both $w_1=0$ and $w_2=0$, which is a contradiction. 

Therefore, it remains to prove this claim. Choose some $s_1 > \bar{s}_1$ 
to be fixed later. We now use the variables $(x,y)$ and set 
\[
r_1 \equiv  \e \cosh (\bar{s}_\e- s_1).
\]
Then the boundary of $[-\bar{s}_\e + s_1, \bar{s}_\e - s_1] \times 
S^1$ consists of two circles of radius $r_1$, one in each of
the regions $II_1$ and $II_2$, which we denote by $\del B_{r_1}^i$. 
Also, set for $i=1,2$
\[
\tilde{u}_{i,\e} \equiv u_i + h_{\e,i}.
\]

Now multiply ${\LL}_\e w_\e =0$ by any one of the `low eigenmode'
Jacobi fields $\Phi = \Phi_\e^{j,\pm }$, $j = 0, \pm 1$, and integrate over 
$[-\bar{s}_\e + s_1, \bar{s}_\e - s_1] \times S^1$. If we set
\begin{equation}
J^i \equiv  \left( \int_{\del B_{r_1}^i} \frac{\Phi \, \del_r w - w \, \del_r 
\Phi}{ (1+|\nabla \tilde{u}_{i,\e} |^2)^{1/2}} \,  r \, d \theta - \int_{\del 
B_{r_1}^i} \frac{\nabla \tilde{u}_{i,\e} \cdot ( \Phi \, \nabla w - w \, 
\nabla \Phi)}{(1+|\nabla \tilde{u}_{i,\e} |^2)^{3/2}} \del_r \tilde{u}_{i,\e} 
\, r \, d\theta  \right),
\label{eq:ident}
\end{equation}
then we obtain by integration by parts and (\ref{eq:1.5}) that
\begin{equation}
J^1 + J^2 = 0.
\label{eq:sum}
\end{equation}

We substitute in each of the Jacobi fields in turn into this
equality to get different information. First of all, we note
that the estimates for $h_{\e, i}$ in Lemma~\ref{le:7.7} show
that $|\nabla \tilde{u}_{i,\e}| \leq c\, |\nabla u| \leq c'\,
r$. It may then be checked that the first integral in $J^i$
always contains the dominant terms of the expansion with
respect to $r$, and furthermore, that the denominator
$(1+|\nabla \tilde{u}_{i,\e}|^2)^{1/2}$ in this integral 
may be replaced by $1$ without affecting the first two
terms of the expansion.  Hence we shall really be only 
computing the leading asymptotic terms in (\ref{eq:sum})
as $\e \to 0$. Finally, we note that
\[
w_i = (-1)^{i+1} \, \left(\frac{c_1^i \cos \theta + c_2^i\sin\theta}{r}
+ c_0^i\,\log r + {\cal O}(1)\right).
\]

First we set $\Phi = \Phi^{0,+}_\e$. Since $\Phi^{0,+}_\e =
(-1)^i$ in $II^i$, we get that
\[
J^i \sim - \,  \int \left( \frac{-(c_1^i \cos \theta + c_2^i\sin \theta)}{r^2}
+ \frac{c_0^i}{r} + {\cal O}(1) \right) \, r \, d\theta + {\cal O}(r_1).
\]
The coefficient of $r^{-2}$ integrates to zero, and so 
\[
c_0^1  +  c_0^2  =  {\cal O} (e^{- s_1}).
\]
Since, this holds for every $s_1$, we conclude $c_0^1+c_0^2=0$. 

Next let $\Phi =\Phi^{\pm 1,-}$. Using Proposition~\ref{pr:14}, 
the leading terms of the expansion is
\[
J^i \sim 2 \, \int ( c_1^1 \cos \theta + c_2^1 \sin\theta )S \, d\theta +{\cal O}( r_1),
\]
where $S$ is equal to either $\cos \theta$ or $\sin \theta$. This
gives 
\[
c^1_1 + c^2_1=c^1_2 + c^2_2=0.
\]
When $\Phi = \Phi^{\pm 1,+}$ then we use Proposition~\ref{pr:4.11}
along with the fact that $u_i$ may be approximate by
its second order Taylor polynomial and (as before) $h_{i,\e}$
may be disregarded. This gives
\[
c_1^1 \, \del^2_{xx}u_1 + c_2^1 \, \del_{xy}^2 u_1 + c_1^2 \, 
\del^2_{xx}u_2 + c_2^2 \, \del_{xy}^2 u_2 =0, \qquad
\mbox{\rm and } \qquad
\]
\[ 
c_1^1 \, \del^2_{xy}u_1 + c_2^1 \, \del_{yy}^2 u_1 + c_1^2 \, \del^2_{xy}u_2 
+ c_2^2 \, \del_{yy}^2 u_2=0,
\]
with all partial derivatives computed at the origin. We write these 
equations all together as
\[
\left(
\begin{array}{ccccccc}
1  & 0  &  1  & 0  \\[3mm]
0  & 1  &  0   & 1 \\[3mm]
\del_{xx}^2u_1 & \del_{xy}^2u_1  & \del_{xx}^2 u_2 &  \del_{xy}^2u_2   
\\[3mm]
\del_{xy}^2u_1 & \del_{yy}^2u_1  & \del_{xy}^2 u_2 &  \del_{yy}^2u_2   
\end{array}
\right)
\left(
\begin{array}{rlllllll}
c^1_1 \\[3mm]
c^1_2 \\[3mm]
c^2_1 \\[3mm]
c^2_2
\end{array}
\right) =
\left(
\begin{array}{rlllllll}
0 \\[3mm]
0 \\[3mm]
0 \\[3mm]
0
\end{array}
\right)
\]
Since we are assuming that (\ref{eq:nondegeqn}) does not hold, 
this matrix is not singular and so $c^1_1=c^1_2=c^2_1=c^2_2=0$.

Finally, we let $\Phi = \Phi_\e^{0,-}$. We have already shown that
$c_j^i = 0$ when $j = 1,2$, and the leading term of 
$\Phi$ is $\log \e$. Then the leading singular term in the expansion
for $J^i$ is 
\[
\ds \frac{J^i}{\log \e} \sim  (-1)^{i+1} \, \int c_0^i \, d \theta + {\cal O} ( r_1).
\]
Hence we get $c_0^1=c_0^2$, which together with the
fact that $c_0^1+c_0^2=0$, implies that $c_0^1=c_0^2=0$.
The claim, and the theorem, is now proved. \hfill $\Box$

\ms

Using this result, the proof of Theorem~\ref{th:last} is now easy 
to complete. In fact, we merely translate (\ref{eq:nondegeqn})
into a more explicit equation involving the principle curvatures
of the surfaces $\Si_i$ at the points $p_i$ and the angle
$\theta$. We shall denote the principle curvatures of $\Si_i$
by $\alpha_i$ and $\beta_i$.

Recall that we had oriented the surfaces so that the $x$ and $y$
axes are principle directions for $\Si_1$. Thus
\[
\del_{xx}^2 u_1(0,0) = \alpha_1, \qquad \del_{xy}^2 u_1(0,0) = 0,
\qquad \mbox{ and } \quad \del_{yy}^2 u_1(0,0) = \beta_1.
\]
On the other hand, using coordinates $(\tilde{x},\tilde{y})$ 
defined by $\tilde{x} = \cos \theta \, x + \sin \theta  \, y$ and
$\tilde{y} = -\sin \theta  \, x + \cos \theta \, y$, we conclude that
\[
\del_{xx}^2 u_2(0,0) = -\alpha_2 \cos^2 \theta - \beta_2 \sin^2 \theta, 
\qquad \del_{xy}^2 u_2(0,0) = (\beta_2 - \alpha_2) \sin \theta \cos\theta,
\]
\[
\qquad \mbox{ and } \quad \del_{yy}^2 u_2(0,0) = -\alpha_2 \sin^2 \theta
-\beta_2 \cos^2 \theta.
\]
(Recall that $\Si_2$ is oppositely oriented to $\Si_1$, which
accounts for the change of signs.) 

We have now proved that the surface $S_\e$ can be degenerate
for $\e$ sufficiently small only if 
\[
(\alpha_2 \cos^2 \theta + \beta_2 \sin^2 \theta + \alpha_1)
(\alpha_2 \sin^2 \theta + \beta_2 \cos^2 \theta + \beta_1)
- \sin^2 \theta \cos^2 \theta (\beta_2 - \alpha_2)^2 = 0.
\]
Some algebra shows that this is equivalent to
\[
(\alpha_1 \alpha_2 + \beta_1 \beta_2) \sin^2 \theta + 
(\alpha_1 \beta_2 + \beta_1 \alpha_2) \cos^2 \theta
+ (\alpha_1 \beta_1 + \alpha_2 \beta_2) = 0.
\]
This is equivalent to a quadratic polynomial in $\cos \theta$, and
hence either the polynomial is identically satisfied, or
else there are at most two values of $\cos \theta$ for which
it  vanishes (and hence at most four values of $\theta$).
The polynomial can only be identically satisfied if
\[
\alpha_1 \alpha_2 + \beta_1 \beta_2 = 
\alpha_1 \beta_2 + \beta_1 \alpha_2 =
- (\alpha_1 \beta_1 + \alpha_2 \beta_2).
\]
Recalling that $\beta_i = H_0 - \alpha_i$, the first equality 
implies that $\alpha_1 H_0 + \alpha_2 H_0= 2 \alpha_1 \alpha_2 + H_0^2/2$
while the second gives $(\alpha_1 + \alpha_2)(2H_0 - \alpha_1 -
\alpha_2) = 0$. These equations together yield 
$\alpha_1 = H_0/2$, $\alpha_2 = 3H_0/2$, and hence $\beta_1 = H_0/2$,
$\beta_2 = -H_0/2$, so that $\Si_1$ is umbilic at $p_1$,
or else $\Si_2$ is umbilic at $p_2$ (with principal curvatures
$(H_0/2,H_0/2)$) while the principal curvatures of $\Si_1$
at $p_1$ are $3H_0/2$ and $-H_0/2$. 

To proceed further, we show that the set of points in
$\Si_1 \times \Si_2$ where the principal curvatures
can have these set values is no more than three dimensional.
The real analyticity of CMC surfaces shows that the locus
of points with fixed principal curvatures is an analytic set,
hence either a discrete set, a collection of analytic
arcs or else the whole surface. Now by definition an
{\it isoparametric surface} $\Si$ is one for which the
principal curvatures are everywhere constant. It is a classical 
theorem of Cartan that the only isoparametric surfaces
in ${\R}^3$, even locally, are subdomains of the sphere
and the cylinder. Hence although $\Si_1$ could be everywhere
umbilic, it is impossible for $\Si_2$ to have principal
curvatures $(3H_0/2,-H_0/2)$ on an open set. This shows
that the portion of the degeneracy set ${\cal S}$ which
includes the complete $S^1$ factor lies over a set in
$\Si_1 \times \Si_2$ which is three dimensional if one of the 
$\Si_i$ is a subdomain of the sphere, and at most two dimensional
otherwise. 

The proof of Proposition~\ref{pr:2} is now complete. 

\section{Appendix I : Using different vector fields to parameterize 
all nearby surfaces}

This section is entirely taken from \cite{MPa}. We have included it here for 
the sake of completeness. Let $\Sigma$ be a regular orientable surface, with 
unit normal vector field $N$. Suppose that $\bar{N}$ is another unit vector 
field along $\Sigma$ which is nowhere tangential. By the inverse function 
theorem, for any $p_0 \in S$ there are neighbourhoods 
${\cal U}$ and ${\cal V}$ near $(p_0,0)$ in $\Sigma \times {\R}$ and 
a diffeomorphism $(\phi(p,s),\psi(p,s))$ from ${\cal U}$ to ${\cal V}$ 
such that 
\begin{equation}
p + s N(p) = \phi(p,s) + \psi(p,s)\bar{N}(\phi(p,s)).
\label{eq:1.88}
\end{equation}
Here $\phi(p,0) = p$ and $\psi(p,0) = 0$. To determine the first order Taylor 
series of these functions in $s$, differentiate (\ref{eq:1.88}) with 
respect to $s$ and set $s=0$. This gives 
\[
N (p) = \frac{\del \phi}{\del s}(p,0) + 
\frac{\del \psi}{\del s}(p,0)\bar{N}(p),
\]
and so, taking the normal component of this, we get
\[
1 = \frac{\del \psi}{\del s}(p,0)\, N (p)\cdot \bar{N}(p), 
\qquad \mbox{\rm or } 
\quad \frac{\del \psi}{\del s}(p,0) = 1/\big( N(p)\cdot \bar{N}(p)\big).
\]
Hence
\[
\psi(p,s) = \frac{s}{N(p)\cdot \bar{N}(p)} + O(s^2).
\]
On the other hand, taking the tangential component and using this expansion
of $\psi$ yields
\[
0 = \frac{\del \phi}{\del s}(p,0) + \frac{s}{N(p)\cdot \bar{N}(p)}
\bar{N}_t(p),
\]
where $\bar{N}_t(p)$ is the tangential component of $\bar{N}$. Thus
\[
\phi(p,s) = p - \frac{s}{N (p)\cdot \bar{N}(p)} \bar{N}_t(p) + O(s^2).
\]

Next, any  ${\cal C}^2$ surface close to $\Sigma$ can be parameterized either 
as a normal graph of some function $w$ over $S$, using the vector field $N$, 
or as a graph of a different function $\bar{w}$ using the vector field 
$\bar{N}$. These functions are related by
\[
p + w(p) \, N(p) = \bar{p} + \bar{w}(\bar{p}) \, \bar{N}(\bar{p}) = 
\phi(p,w(p)) + \psi(p,w(p)) \, \bar{N}(\phi(p,w(p))).
\]
Using the expansions above, we see that $\bar{w}(\bar{p}) = w(p)/ (N(p) \cdot 
\bar{N}(p)) + O(\|w \|^2)$. 

The mean curvature operators on these two functions, which we call 
$H_{w}$ and $\bar{H}_{\bar{w}}$, respectively, are related by 
\begin{equation}
\bar{H}_{\bar{w}}(\bar{p}) = H_{w}(p).
\label{eq:1.9}
\end{equation}
Differentiating this with respect to $\bar{w}$ and setting $\bar{w} =
0$, we get
\begin{equation}
D_{\bar{w}} \bar{H}_{0}(u) = D_{w} H_{0}(\bar{N}\cdot N \, u) + 
\left(\nabla H_{0} \cdot \bar{N}_t\right) \, u,
\label{eq:1.10}
\end{equation}
for any scalar function $u$. In the special case where the surface 
$\Sigma$ has constant mean curvature, this reduces to 
\begin{equation}
D_{\bar{w}} \bar{H}_{0} (u) = D_{w} H_{0}( \bar{N} \cdot  N \, u ) .
\label{eq:1.11}
\end{equation}

\section{Appendix II : Precise expansions for the mean curvature operator}

\subsection{Proof of the expansion (\ref{eq:1.8})}

We use polar coordinates $(t, \theta)$ and write
$u(t, \theta)= u(e^{t} \, \cos \theta, e^{t} \, \sin \theta)$.
With this notation, (\ref{eq:1.2}) becomes $\nabla^k u (t, \theta) =
{\cal O} (e^{-2t})$, where $\nabla$ is now the gradient with respect to
$t$ and $\theta$. In the same way, the fact that the function $w$ as well
as all its derivatives are assumed to be small  means that $\nabla^k w 
(t, \theta) = {\cal O} (e^{-t})$.

The mean curvature operator is
\[
H_{u+w} = e^{t} \nabla \ds\left(  \frac{ e^{t} \nabla (u+w)}{\left( 
1+ e^{2t} |\nabla u+w|^2\right)^{1/2} } \right) .
\] 
We start by noting that
\[
\ds \frac{1}{( 1+ e^{2t} |\nabla u+w|^2)^{1/2}}  =   \frac{1}{( 1+ e^{2t} 
|\nabla u|^2 )^{1/2}} - e^{2t} \ds \frac{\nabla u \cdot \nabla w}{
( 1+ e^{2t} |\nabla u |^2 )^{3/2}} + Q_1 (e^t \, \nabla w ) ,
\]
where the function $Q_1(\cdot)$ is a function all of whose derivatives 
are bounded in ${\cal C}^k ([-\log \rho, +\infty) \times
S^1)$ and which satisfies $Q_1 (0)=0$ and $\nabla Q_1(0)=0$.
 
Next, $H_{u+w}$ is given by
\[
\ds  e^{t} \nabla \left( \frac{e^t}{( 1 + e^{2t} 
|\nabla u|^2 )^{1/2}} \, \nabla (u+w) - e^{3t} \, \nabla (u+w) \, 
\ds \frac{\nabla u \cdot \nabla w}{( 1+ e^{2t} |\nabla u |^2)^{3/2}} + 
e^{t} \nabla (u+w) \,  Q_1 \left( e^t \, \nabla w
\right)\right) .
\]
From this it follows at once that
\[
H_{u+w} = H_u + \Lambda_u w -Q'_u (e^t \, \nabla w, e^t \, \nabla^2 w) -
e^t \, Q''_u(e^t \, \nabla w, e^t \, \nabla^2 w) ,
\]
where $Q'_u(\cdot, \cdot)$ and $Q''_u(\cdot, \cdot)$ satisfy the required
properties as stated in (\ref{eq:1.8}). 

\subsection{Proof of the expansion (11)}

We consider a surface parameterized by
\[
(s, \theta) \longrightarrow \bx_\e^c (s, \theta) + \e \, \cosh s \, 
\tilde{w}(s, \theta) \, n (s, \theta),
\]
for some regular function $w(s, \theta) =\e \, \cosh s \, \tilde{w}(s, 
\theta)$, where $\bx_\e^c$ is given by (\ref{eq:3.3}) and $n (s, \theta)$
is given by (\ref{eq:5.2}). 
A simple computation shows that the 
coefficients of the first fundamental form are then given by
\[
E_w = \e^2 \, \left( \cosh^2 s - 2 \cosh s \, \tilde{w}  + 
\cosh^2 s \, (\tilde{w}^2 +\tilde{w}_s^2)
+ 2 \sinh s \, \cosh s  \, \tilde{w} \, \tilde{w}_s \right), 
\]
\[
F_w =  \e^2 \left( \sinh s \, \cosh s \, \tilde{w}_\theta \, 
\tilde{w} + \cosh^2 s \, \tilde{w}_s
\, \tilde{w}_\theta \right)
\]
and
\[
G_w = \e^2 \, \left( \cosh^2 s + 2 \cosh s \, \tilde{w}  + 
\tilde{w}^2 + \cosh^2 s \, 
\tilde{w}_\theta^2 \right). 
\]
Notice that these can be written  as
\[
E_w = \e^2 \, \left( \cosh^2 s - 2 \cosh s \, \tilde{w} + \cosh^2 s \, 
P_E ( \tilde{w}, \nabla \tilde{w})
\right)
\]
\[
F_w = \e^2 \, \cosh^2 s \, P_F ( \tilde{w}, \nabla \tilde{w})
\]
\[
G_w = \e^2 \, \left( \cosh^2 s + 2 \cosh s \, \tilde{w} + \cosh^2 s \, 
P_G ( \tilde{w}, \nabla \tilde{w})
 \right)
\]
where $P_E, P_F$ and $P_G$ are polynomials, homogeneous of
degree $2$, whose coefficients are bounded functions of $s$ and 
$\theta$. 
Moreover, using this we may write
\[
E_w G_w -F_w^2 = \e^4 \, \cosh^4 s \,  \left( 1  + P_{EG-F^2} 
(\tilde{w}, \nabla \tilde{w})
 \right),
\]
Where $P_{EG-F^2}$ is a polynomial consisting of terms homogeneous of degree 
2 or 4, whose coefficients are bounded functions of $s$ and 
$\theta$. 

In the same way, we compute the coefficients 
of the second fundamental form and find that these are given by 
\[
\sqrt{E_w G_w -F_w^2} \, e_w = -\e^3 \, \left(
\cosh^2 s + \cosh^3 s \, \tilde{w}_{ss} + \cosh^2 s \, \sinh s \, \tilde{w}_s
+ \cosh^2 s \, P_e (\tilde{w}, \nabla \tilde{w},  \nabla^2 \tilde{w})
\right),
\]
\[
\sqrt{E_w G_w -F_w^2} \,  f_w = - \e^3 \,  \left(\cosh^3 s \, \tilde{w}_{s\theta} +
+ \cosh^2 s \, P_f (\tilde{w}, \nabla \tilde{w},  \nabla^2 \tilde{w})
\right)
\]
and 
\[
\begin{array}{rlll}
\sqrt{E_w G_w -F_w^2} \,  g_w  = - \e^3  & \left(- \cosh^2 s + \cosh^3 s \, 
\tilde{w}_{\theta\theta} + (\cosh^3 s - 2\cosh s) \, \tilde{w} +
\cosh^2 s \, \sinh s \, \tilde{w}_s \right. \\[3mm]
& \quad \left. + \cosh^2 \, P_g (\tilde{w}, \nabla \tilde{w},  \nabla^2 \tilde{w})
\right),
\end{array}
\]
where, here also, $P_e, P_f, P_g$ are polynomials without any constant nor 
any linear terms and  all of whose coefficients are bounded functions of $s$ 
and $\theta$.
 
The mean curvature operator may then be expressed in terms of these
coefficients as 
\[
H_{w}  = \frac{e_w G_w - 2 f_w F_w + g_w E_w}{E_w G_w -F_w^2}.
\]
Using the previous expansions we obtain
\[
\begin{array}{rlll}
H_w & =   \ds - \frac{1}{\e\, \cosh s} & \ds  \left(  
\tilde{w}_{ss} +\tilde{w}_{\theta \theta} + 2 \tanh s \, \tilde{w}_s + 
\tilde{w} + \frac{2}{\cosh^2 s} \tilde{w} \right. \\[3mm]
& \, & \ds \left. \frac{1}{\cosh s}\, P' (\tilde{w}) +
P''(\tilde{w})\right) \left( 1 + \tilde{P}' (\tilde{w})\right) ,
\end{array}
\] 
where $P'$, $\tilde{P}'$ and $P''$ are functions of $\tilde{w}$,
$\nabla\tilde{w}$ and $\nabla^2\tilde{w}$, all of whose  
partial derivatives are 
bounded functions in ${\cal C}^k ([-s_\e, s_\e] \times S^1)$, for all 
$k\geq 0$, uniformly in $\e$. Moreover, these functions  satisfy
\[
P' (0,0,0)= \tilde{P}'(0,0,0)= P'' (0,0,0) = 0 
\]
\[
\nabla P' (0,0,0)= \nabla \tilde{P}' (0,0,0) = \nabla P'' (0,0,0) =0 
\]
and $P''$ satisfies in addition  
\[
 \nabla^2 P'' (0,0,0)=0.
\]

After having performed the change of function 
$\tilde{w} = w/(\e \cosh s)$, we conclude that
\[
\begin{array}{rlll}
H_w = - \ds \frac{1}{\e^2 \, \cosh^2 s} {\cal L} w & + 
\ds \frac{1}{\e \cosh^2 s} Q_\e'\left( \frac{w}{\e \cosh s},\frac{\nabla w}{\e
\cosh s},\frac{\nabla^2 w}{\e \cosh s}\right) \\[3mm]
& + \ds \frac{1}{\e \cosh s} Q_\e'' \left( \frac{w}{\e \cosh s},\frac{\nabla
w}{\e \cosh s},\frac{\nabla^2 w}{\e \cosh s}\right), 
\end{array}
\]
where $Q'_\e$ and $Q''_\e$ are functions all of whose  partial 
derivatives are bounded in ${\cal C}^k ([-s_\e, s_\e] \times S^1)$, 
for all $k\geq 0$, uniformly in $\e$. 
Moreover, these functions  satisfy
\[
Q'_\e (0,0,0)= Q''_\e (0,0,0) = 0 \qquad \mbox{and}\qquad \nabla 
Q'_\e (0,0,0)= \nabla Q''_\e (0,0,0) =0 .
\]
and $Q''_\e$ satisfies in addition  
\[
 \nabla^2 Q''_\e (0,0,0)=0.
\]

\section{Appendix III : 
Mapping properties of the Laplace operator in a cylinder}

We collect here various results whose proofs are slight 
modifications of the proof of Proposition~\ref{pr:6}.
\begin{lemma}
Assume that $\delta \in (1,2)$ and that $0 < s_0 < s_1$. Then, there exists 
some operator 
\[
G_{s_0, s_1} : {\cal C}^{0,\alpha}_{\delta +2}([s_0 , s_1] \times S^1)
\longrightarrow  {\mathcal C}^{2,\alpha}_{\delta}([s_0 , s_1] \times S^1),
\]
such that, for all $f \in {\cal C}^{0,\alpha}_{\delta+2 }([s_0 , s_1] 
\times S^1)$, the function $w= G_{s_0, s_1} (f)$ is a solution of the problem
\[
\left\{ 
\begin{array}{rllll} 
e^{2s} \Delta w & = & \ds  f \qquad & \mbox{in} \quad (s_0, s_1) \times 
S^1 \\[3mm] 
\pi'' w  & =  &  0 \qquad & \mbox{on}\quad \{s_1\} \times  S^1\\[3mm]
w  & =  &  0 \qquad & \mbox{on}\quad \{s_0\} \times S^1.
\end{array} 
\right. 
\]
In addition, we have $|| G_{s_0, s_1} (f) ||_{2, \alpha, \delta} \leq c \,
|| f ||_{0, \alpha, \delta + 2} ,$ for some constant $c>0$ independent of
$s_0, s_1$.
\label{le:2}
\end{lemma}
{\bf Proof :} Using separation of variables as in the proof of 
Proposition~\ref{pr:6}, we now write
\[
\tilde{w} = \sum_{n \in {\Z}} w_n (s) e^{in \theta} \qquad \mbox{and} 
\qquad  f =  e^{2s}\,\sum_{n \in {\Z}}  f_n (s) e^{in \theta}.
\]
By linearity, we may assume that $ | f| (s) \leq e^{(\delta + 2) s}$ and 
therefore we find  $|f_n |(t) \leq  e^{\delta  s}$.
This time, for all $|n|\geq 2$, we see that $w_n$ has to solve 
\[
\ddot{w}_n -n^2 w_n = f_n \qquad \mbox{in} \qquad (s_0, s_1),
\]
and $w_n (s_0)=w_n(s_1)=0$. It is easy to see that the function $
\ds \frac{2}{n^2 -\delta^2}\, e^{\delta s}$ is a supersolution for our problem 
therefore this yields, for all $|n|\geq 2$
\[
|w_n | (s) \leq \displaystyle \frac{2}{n^2 -\delta^2} \, e^{\delta s}.
\]
For $n=0$ and $n=\pm 1$, we use the explicit formula 
\[
w_{0} (s) = \int_{s_0}^s  \int_{s_0}^t  f_{0} (u) \, du\, dt
\qquad \mbox{and}\qquad 
w_{\pm 1} (s) = e^{-s} \int_{s_0}^s e^{2 t} \int_{s_0}^t e^{- u} \, f_{\pm 1} 
(u) \, du \, dt.
\]
Summation over $n$ and Schauder's estimates lead to the desired result. 
\hfill $\Box$

\ms

Our next Lemma is a variant of the previous result. 
\begin{lemma}
Assume that $\delta \in (0,1)$ and that $0 < s_0 < s_1$. Then, there exists 
some operator 
\[
G'_{s_0, s_1} : {\cal C}^{0,\alpha}_{\delta +2}([s_0 , s_1] \times S^1)
\longrightarrow  {\mathcal C}^{2,\alpha}_{\delta}([s_0 , s_1] \times S^1),
\]
such that, for all $f \in {\cal C}^{0,\alpha}_{\delta+2 }([s_0 , s_1] 
\times S^1)$, the function $w= G'_{s_0, s_1} (f)$ is a solution of the problem
\[
\left\{ 
\begin{array}{rllll} 
e^{2s} \Delta w & = & \ds  f \qquad & \mbox{in} \quad (s_0, s_1) \times 
S^1 \\[3mm] 
\pi' w  & =  &  0 \qquad & \mbox{on}\quad \{s_1\} \times  S^1\\[3mm]
w  & =  &  0 \qquad & \mbox{on}\quad \{s_0\} \times S^1.
\end{array} 
\right. 
\]
In addition, we have $|| G'_{s_0, s_1} (f) ||_{2, \alpha, \delta} \leq c \,
|| f ||_{0, \alpha, \delta + 2} ,$ for some constant $c>0$ independent of
$s_0, s_1$.
\label{le:2.2}
\end{lemma}
{\bf Proof :} The only difference with the proof of the previous 
result is that, this time $\ds \frac{2}{n^2 -\delta^2}\, e^{\delta s}$
is a supersolution for our problem for all $|n|\geq 1$. \hfill $\Box$

\ms

We will also need 
\begin{lemma}
Assume that $\delta' \in (-1,0)$ and that $s_0 >0$. Then, there 
exists an operator 
\[
\hat{G}_{s_0} : {\cal C}^{0,\alpha}_{\delta'+2}([s_0 , + \infty) \times S^1)
\longrightarrow  {\mathcal C}^{2,\alpha}_{\delta'}([s_0 , +\infty) \times S^1)
\oplus  {\R},
\]
such that, for all $f \in {\cal C}^{0,\alpha}_{\delta'+2}
([s_0 ,+\infty) \times S^1)$, the function $w= \hat{G}_{s_0} (f)$ is the
unique solution of the problem
\[
\left\{ 
\begin{array}{rllll} 
e^{2s} \Delta w & = & \ds  f \qquad & \mbox{in} \quad (s_0, +\infty)
 \times S^1 \\[3mm] 
w  & =  &  0 \qquad & \mbox{on}\quad \{s_0\} \times S^1 ,
\end{array} 
\right. 
\]
which belongs to the space $  {\mathcal C}^{2,\alpha}_{\delta'}([s_0 , +\infty)
\times S^1) \oplus  {\R}$. In addition, if we decompose $w(s, \theta) =
v(s, \theta) + c_0 \in {\mathcal C}^{2,\alpha}_{\delta'}([s_0 , +\infty)
\times S^1) \oplus  {\R}$, we have $e^{-\delta' s_0} |c_0| + || v 
||_{2,\alpha, \delta'} \leq c \, || \tilde{f}||_{0, \alpha, \delta'+2 } $,
for some constant $c>0$ independent of $s_0$.
\label{le:3}
\end{lemma}
{\bf Proof :}
Again, using separation of variables as  
in the proof of Lemma~\ref{le:2}, 
we  write
\[
w = \sum_{n \in {\Z}} w_n (s) e^{in \theta} \qquad \mbox{and} 
\qquad  f = e^{2s} \sum_{n \in {\Z}}  f_n (s) e^{in \theta}.
\]
By linearity, we may assume that $ | f| (s) \leq e^{(\delta'+2) s}$ and 
therefore we find  $|f_n |(t) \leq  e^{\delta'  s}$.
Here, for all $|n|\geq 1$, we see that $w_n$ has to solve 
\[
\ddot{w}_n -n^2 w_n = f_n \qquad \mbox{in} \qquad (s_0, +\infty),
\]
and $w_n (s_0)=0$. It is easy to see that the function $\ds \frac{2}{n^2 
- (\delta')^2} \, e^{\delta' s}$ is a supersolution for our problem 
therefore this yields, for all $|n|\geq 1$
\[
|w_n | (s) \leq \displaystyle \frac{2}{n^2  -(\delta')^2}\, e^{\delta' s}.
\]
For $n=0$, the variation of the constant formula provides us 
with the explicit formula
\[
w_{0} (s) = c_0 +  \int_s^{+\infty}  \int_t^{+\infty}  f_{0} (u) \, du\, dt,
\qquad \mbox{with} \qquad 
c_0 =-  \int_{s_0}^{+\infty}  \int_t^{+\infty}  f_{0} (u) \, du\, dt.
\]
And the desired estimates follow at once by summation over $n$ and
direct estimate for $c_0$.\hfill $\Box$

\ms

Using similar arguments, we can also prove  
\begin{lemma} 
There exists an operator
\[
{\cal P} : \pi'' \left( {\cal C}^{2,\alpha}(S^1) \right)
\longrightarrow  {\mathcal C}^{2,\alpha}_{-2}((-\infty,0] \times S^1),
\]
such that, for all $\phi'' \in \pi'' \left( {\cal C}^{2,\alpha}(S^1)\right)$, 
the function  $w={\cal P} (\phi'')$ is the unique solution of
\begin{equation} 
\left\{ 
\begin{array}{rllll} 
\Delta w & = & 0 \qquad & \mbox{in}\quad (-\infty , 0) \times S^1\\[3mm] 
 w & =  &  \phi''   \qquad  & \mbox{on}\quad \{ 0 \} \times S^1 ,
\end{array}  
\right. 
\label{eq:5.9}
\end{equation} 
which belongs to the space ${\cal C}^{2,\alpha}_{-2} ((-\infty, 0] \times 
S^1)$. In addition, we have $ || {\cal P}(\phi'') ||_{2, \alpha, - 2} 
\leq c \, ||\phi'' ||_{2,\alpha}$, for some constant $c>0$.
\label{pr:7}
\end{lemma}
{\bf Proof :} As in the proof of the previous Lemma, we decompose 
$\phi''$ into Fourier series $ \phi'' = \sum_{|n|\geq 2} \phi_n 
e^{in \theta}$. The solution $w$ is then explicitly given by $ w = 
\sum_{|n|\geq 2} \phi_n e^{|n| s} e^{in \theta}$, from which it immediately 
follows that 
\[
 |w| (s) \leq 2 \,  e^{2 s} (1 +   \sum_{n\geq 3} e^{(n-2) s}) 
\, || \phi ||_{2, \alpha}.
\]
Therefore, we already obtain $\sup_{s\leq {-1}} 
e^{-2s} |w| (s) \leq c \, || \phi'' ||_{2, \alpha}$. It also follows from 
the explicit formula for $w$ that $\| w(-1, \cdot) \|_{{\cal C}^{2, \alpha}} 
\leq c \, || \phi'' ||_{2, \alpha}$. Using this last estimate as well as the 
fact that $\Delta w = 0$ in  $(-1,0)\times S^1$, we find from Schauder's 
estimates that $\sup_{s\in (-1,0)} e^{-2s} |w| (s) \leq c \, || \phi'' 
||_{2, \alpha}$. The other estimates, for the derivatives of $w$, follow 
again from Schauder's estimates. \hfill $\Box$

\end{document}